\documentclass[review]{elsarticle}

\usepackage{lineno,hyperref}
\modulolinenumbers[5]

\newtheorem{theorem}{Theorem}
\newtheorem{proposition}{Proposition}
\newdefinition{rmk}{Remark}
\newproof{pf}{Proof}
\newproof{pot}{Proof of Theorem \ref{thm2}}

\usepackage{graphicx}
\usepackage{amsmath}
\usepackage{amsfonts}
\usepackage{algorithm}
\usepackage{algorithmic}
\usepackage{multirow}
\usepackage{tabularx}
\usepackage{caption}
\usepackage{subcaption}
\usepackage{tablefootnote}

\usepackage{soul,color}
\usepackage{url}
\soulregister\cite7
\soulregister\ref7

\usepackage{geometry}
\allowdisplaybreaks

\begin{document}

\begin{frontmatter}

\title{Park-and-Ride Facility Location Selection under Nested Logit Demand Function}

\author{Sang Hyun Kim\fnref{myfootnote}}
\address{Korea Aerospace University, Goyang-si, Gyeonggi-do, 10540, Republic of Korea}
\fntext[myfootnote]{Assistant Professor, School of Air Transport, Transportation, and Logistics}

\author{Sangho Shim\fnref{myfootnote2}\corref{mycorrespondingauthor}}
\address{Robert Morris University, Moon Township, PA 15108, USA}
\fntext[myfootnote2]{Associate Professor, School of Engineering, Mathematics and Science}
\cortext[mycorrespondingauthor]{Corresponding author}
\ead{shim@rmu.edu}

\begin{abstract}
Park-and-ride (P\&R) facilities are car parks where users can transfer to public transportation, by which they reach their final destinations. Commuters can use P\&R facilities or choose to travel by car to their destinations, and individual choice behavior is assumed to follow a logit model. The P\&R facility location problem identifies locations for a fixed number of P\&R facilities from among potential locations such that the number of users of the P\&R facilities is maximized. 
This problem has previously been formalized under a multinomial logit (MNL) demand function. However, as it imposes the strong condition of the independence of irrelevant alternatives (IIA), the MNL model is unable to represent the real-world P\&R facility location problem exactly. Respecting the nested structure of individual choice behavior, we generalize the MNL model to a nested logit (NL) model and develop two computational methods---neighborhood search and randomized rounding---to solve large-scale P\&R facility location problems under the NL demand function. The neighborhood search method first finds one feasible solution and then improves the feasible solution to the next one along an edge of the polyhedron whose vertices are the feasible solutions. The neighborhood search method is randomized to create an adaptive randomized rounding procedure. Computational experiments verified that the computational methods were able to solve the nonlinear optimization problem under the NL demand function on 1,000 medium-scale instances to exact optimality rapidly. Specifically, the methods were verified to solve the MNL model to exact optimality 10,000 times faster than the mixed integer linear programming formulation described in the literature. We performed additional computational experiments to assess the performance of our computational methods on a variety of large-scale instances. Our computational analysis also elucidates the difference between the MNL and NL models.
\end{abstract}

\begin{keyword}
traffic\sep park-and-ride facility location problem\sep decision-dependent demand\sep nested logit model\sep randomized rounding
\end{keyword}

\end{frontmatter}


\section{Introduction}\label{sec:intro}
A transit facility is a place that provides access to transit services \citep{law}. This paper focuses on park-and-ride (P\&R) facilities, which include bus, train, and air mobility stations; generally, people drive to reach these facilities, park their cars, and then transfer to public transportation. As the cost of construction for transit facilities is high and it is difficult to change their locations once they are built, the location and number of new transit facilities must be carefully determined through quantitative and qualitative analyses of relevant factors such as cost, demand, and level of service. 
The problem examined in this paper assumes that a large number of new transit facilities are to be built because of the introduction of a new mode of transportation such as urban air mobility. This park-and-ride facility location problem (P\&R FLP) finds a fixed number of optimal P\&R locations, maximizing the total transportation demand served by the P\&R facilities ({\it e.g.}, the total number of P\&R users or the modal share of P\&R facilities).

The P\&R FLP is a $p$-hub location problem whose objective is to identify locations for a fixed number ($p$) of facilities (in this paper, $p$ is denoted by $N$ to avoid confusion with a probability $p$). \citet{goldman1969p} first addressed the network hub location problem and  \citet{okelly1987p} presented the first 
mathematical formulation of a hub location problem through a study of airline passenger networks. Later, \citet{campbell1996p} presented a mixed integer linear programming (MILP) formulation of the $p$-hub median problem. Based on the $p$-hub formulation, \citet{aros2013p} proposed a $p$-hub approach that utilizes a spatial optimization model taking into account origin--destination (O-D) trips 
and considers transportation demand to determine the proportion of users patronizing the facility. 
{In this study, we solve the P\&R FLP introduced by \citet{aros2013p} by generalizing the objective function of the multinomial logit demand to the nested logit demand function.}

\begin{figure}
    \centering
    \includegraphics[width=0.7\textwidth]{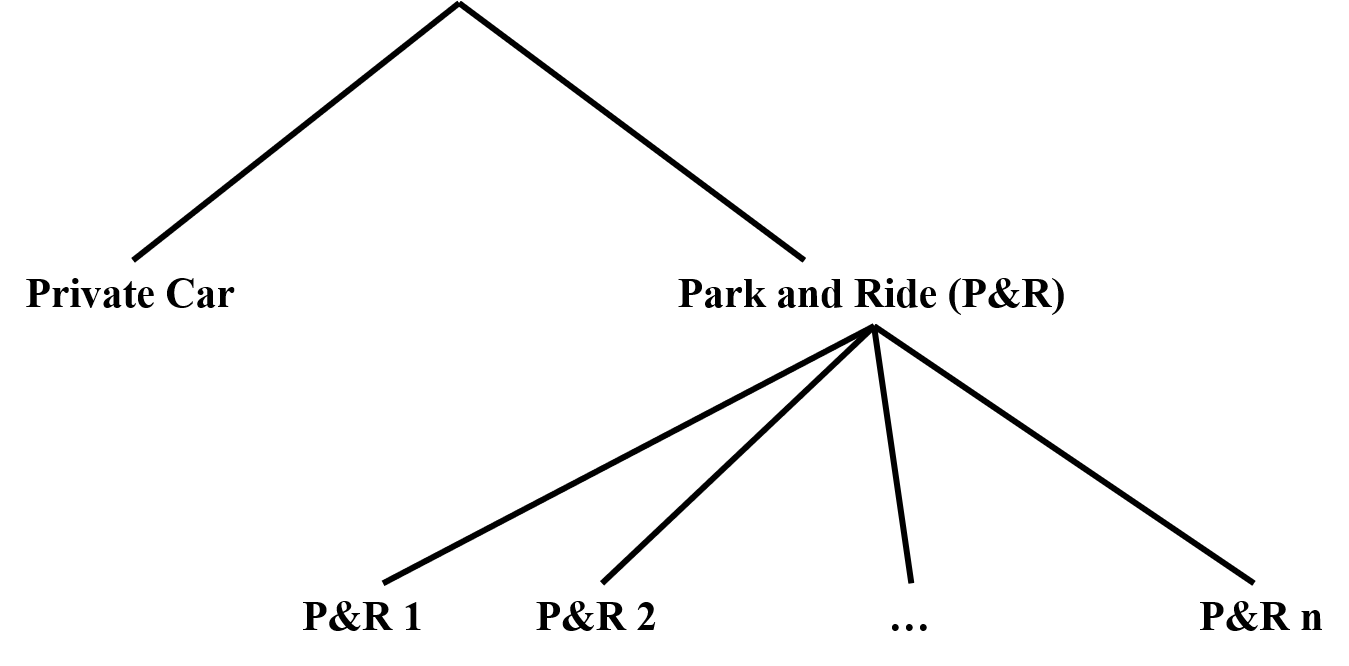}
    \caption{Nesting Structure of Private Car and Park \& Ride}
    \label{f:nest}
\end{figure}

Typically, transportation demand is 
analyzed using a four-step model that comprises trip generation, trip distribution, mode choice, and route assignment \citep{mcnally2007four}. {The first step generates the total number of origin and destination trips, and the second step determines the O-D trips by matching origins with destinations. (For this study, we assumed that the O-D trips are publicly known from previous studies such as those using the Korea Transport Database \citep{ktdb} or from other published works.)} Then, given the O-D trips, the modal share of a particular transportation mode ({\it i.e.}, the probability that that mode is chosen) is analyzed, usually using a multinomial logit (MNL) model because of its flexibility and computational efficiency \citep{mcfadden1977application}. {A distinguishing property of the MNL model is the independence of irrelevant alternatives (IIA): the ratio of the probabilities that two respective alternatives are chosen is independent of the existence or characteristics of other alternatives \citep{mcfadden1974conditional}.} However, the IIA assumption is invalid if some alternatives share unobservable attributes \citep{ben1985discrete}. As a result, the probabilities for choosing such alternatives (those sharing some (unmodeled) attributes) are overestimated by the MNL model. The well-known red/blue bus problem provides an example of such overestimation of probabilities when buses are being chosen \citep{mcfadden1974conditional}.

\begin{figure}
\centering
\begin{subfigure}{.45\textwidth}
    \centering
    \includegraphics[width=\linewidth]{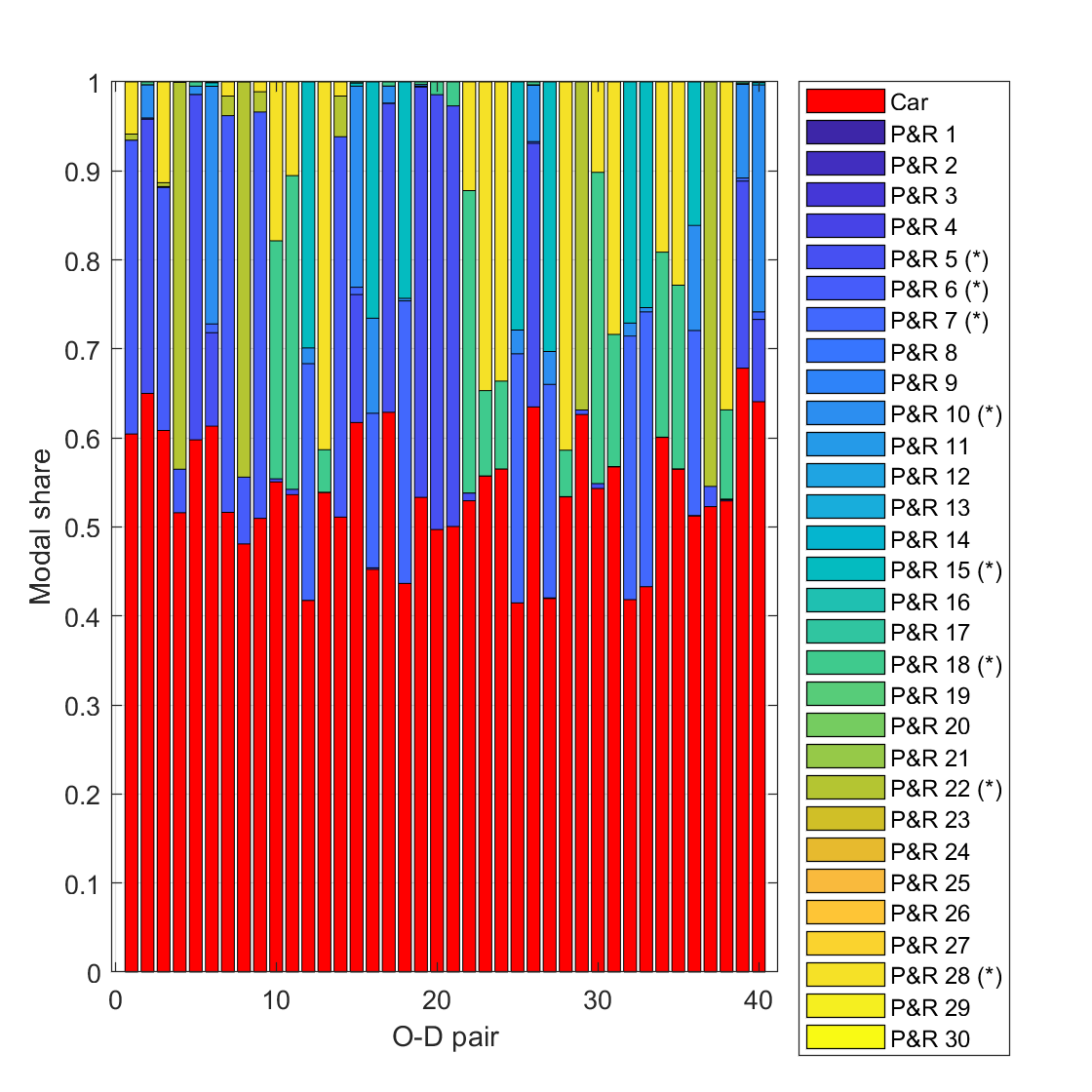}
    \caption{Nested Logit Model}
    \label{f:nl}
\end{subfigure}
\begin{subfigure}{.45\textwidth}
    \centering
    \includegraphics[width=\linewidth]{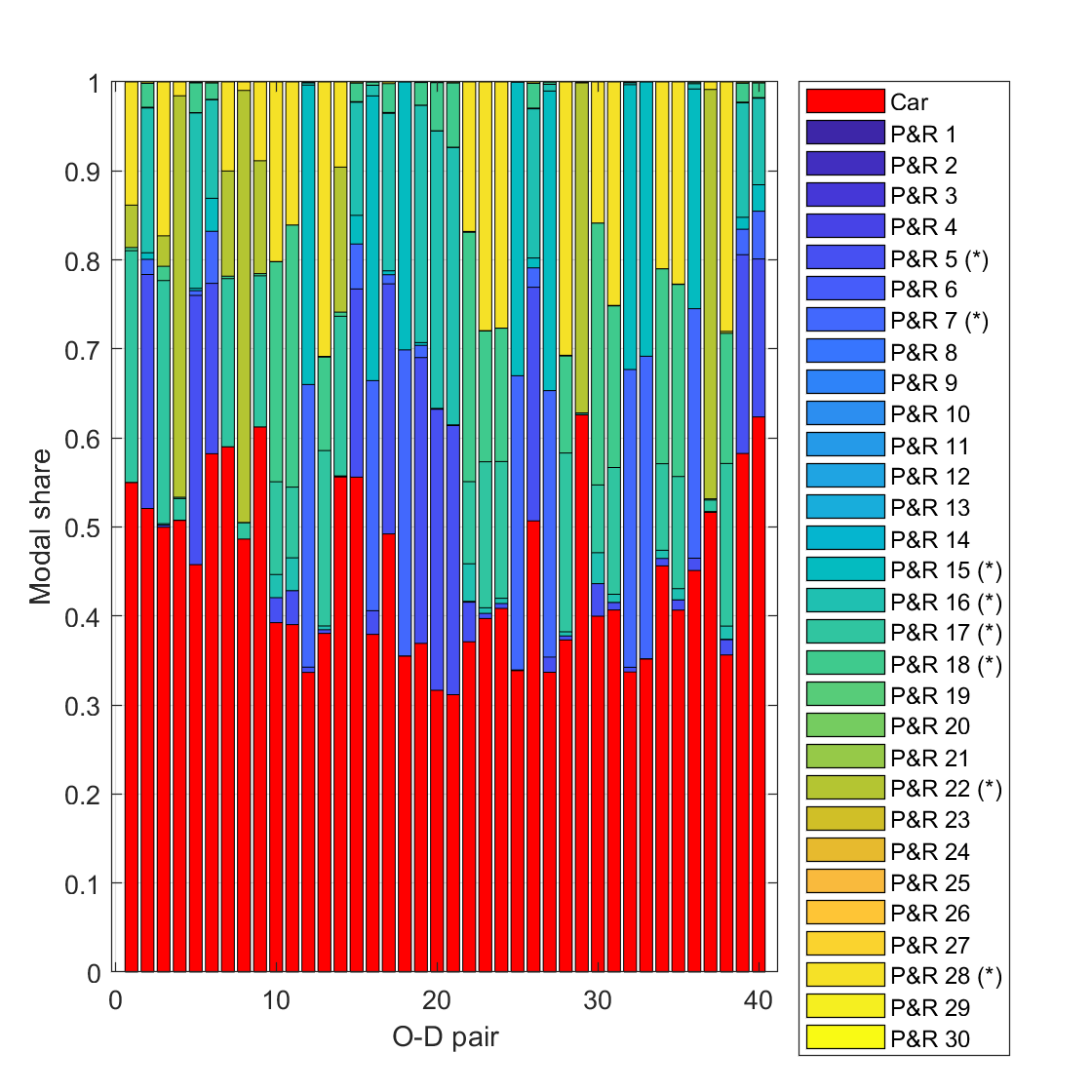}
    \caption{Multinomial Logit Model}
    \label{f:mnl}
\end{subfigure}
\caption{Modal Share and Optimal P\&R Locations (asterisk) of an Instance of Experiment 1}
\label{f:optimal}
\end{figure}

To resolve the IIA property issue, the nested logit (NL) model was proposed, which reflects correlations between alternatives \citep{ben1973structure, mcfadden1978modeling}. 
The NL model is 
represented by a tree structure that groups alternatives that have unobserved attributes in common.
Fig.~\ref{f:nest} shows a nested structure of private car and P\&R options, where the P\&R facilities are grouped in the same subtree. 
When there exists such a nested structure of transit options based on common unobserved attributes, the NL model estimates modal shares better than the MNL model. 
Fig.~\ref{f:optimal} shows that the MNL model underestimates the private car demand, which is indicated by the red bars.\footnote{Fig.~\ref{f:optimal} compares the modal shares estimated by NL and MNL models for 40 commuters that have a private car option plus the eight transit options marked by asterisks. The locations of these eight P\&R options have been optimally selected from among 30 candidates. (See Experiment~1 in Section~\ref{sec:experiments}.)}
A commuter first chooses between using a private car or public transportation; the commuter who decides to use public transportation then chooses one of the P\&R facilities.
Ignoring this nested structure of individual choice behavior, the MNL model treats all nine options (the private car and the eight P\&R options) as being on the same level. 

Together with certain constraints, logit models can be used for the objective function of FLPs. \citet{benati2002maximum} used an MNL model to formulate the objective function of a competitive FLP and computed its upper bounds. \citet{aros2013p} utilized another MNL model to calculate decision-dependent demand for P\&R facilities, which constitutes their objective function. They developed MILP formulations that assume the IIA property for MNL models. \citet{jokar2017facility} modified the linearized model of \citet{aros2013p} to formulate a facility-and-hub location problem. \citet{lopez2021design} estimated customer demand with a constrained MNL model and solved an FLP by considering costs and quality of service. 

However, as \citet{haase2013management} pointed out, there are limitations associated with the IIA property of the MNL model, specifically in spatial choice applications. They utilized a mixed MNL model for a free school-choice problem without assuming any particular distribution for random utility. Using the general random utility model, they simulated students based on population characteristics and solved the problem for the simulated students. \citet{basciftci2021distributionally} considered the effect of facility location on demand uncertainty and proposed a distributionally robust optimization model that minimizes worst-case costs. 
Other transportation researchers \cite{yun1997analysis, kim2005model, park2007nested} reached a consensus that the NL model is more appropriate than the MNL model for the passenger travel data from the Korea Transport Database \cite{ktdb}, especially when the dependence between alternatives cannot be ignored.

{
As \citet{lopez2021design} stated, the more accurate the estimation of transportation demand, the better the decision regarding facility locations. Because it is apparent that P\&R facilities share common unobserved attributes, an NL model can estimate transportation demand more accurately than an MNL model. 
With this in mind, we propose an approach for solving the P\&R FLP introduced by \citet{aros2013p} by generalizing the objective function of the MNL demand to the NL demand function.
Based on the geometry of the polyhedron whose vertices are the feasible solutions, this study developed two computational methods, neighborhood search and randomized rounding, to solve the P\&R FLP under the NL demand function.
}

{
Like the simplex method, the neighborhood search method starts at one vertex of the polyhedron whose vertices are the feasible solutions and moves along an edge to the next vertex (called a neighbor), thereby improving the demand function. 
The neighborhood search method iterates the process of moving from one vertex to the next along an edge until there is no better neighbor remaining.
Based on the same polyhedral geometry, the adaptive randomized rounding procedure randomizes the neighborhood search method, searching for a better random solution (an approximate neighbor) near the best solution found so far.
}

{
In this study, computational experiments were performed to verify that the two computational methods solve the nonlinear optimization problem for medium-scale instances under the NL demand function to exact optimality extremely rapidly.
It was also verified that the methods solve P\&R FLPs under the MNL demand function to exact optimality much faster than the MILP formulation of the MNL model introduced by \citet{aros2013p}. 
We performed additional computational experiments to assess the performance of the proposed computational methods on a variety of large-scale instances for which solutions have not been efficiently found in the literature. Our computational experiments also examine the difference between the MNL and the NL models.
}    
    
In Section~\ref{sec:models}, we develop an NL model as a generalization of the MNL demand function of the P\&R FLP introduced by \citet{aros2013p}; then, we discuss the relationship between the MNL model and the NL model and revisit the MILP formulation developed by \citet{aros2013p} for the MNL model. {Section~\ref{s:neighbor} explores the geometry of the polyhedron whose vertices are the feasible solutions, and a neighborhood search method is developed for solving the P\&R FLP under the NL demand function.} In Section~\ref{sec:heuristics}, the neighborhood search method is randomized, and the adaptive randomized rounding procedure is developed. 
{
Section~\ref{sec:experiments} describes the conditions and reports the results of the computational experiments conducted to assess the performance of the proposed methods; it also presents a sensitivity analysis of the correlation between P\&R locations. Finally, Section~\ref{sec:conclusion} summarizes the findings of the study and presents our conclusions.
}

\section{Nested Logit Demand Function}\label{sec:models}
\subsection{Nested Logit Model}
The NL model can estimate modal shares better than the MNL model can when travelers choose a transit option from among multiple transit options that share common unobserved attributes. The NL model is represented by a tree structure, in which alternatives having unobserved attributes in common are grouped. For example, bus and subway options can be grouped together as public transportation, and their P\&R facilities may be grouped in the same subtree. Fig.~\ref{f:nest} shows an example of a nested structure that consists solely of private car and P\&R options. For simplicity of modeling, we assume that each traveler either uses a private car or chooses to park and ride and that no other mode of transportation is available. 
(The proposed model can easily be extended to a general model with multiple transportation modes, such as one consisting of multiple subtrees representing a private car, P\&R facilities, buses, and subways.)
Private cars move travelers directly from origin to destination. In contrast, travelers who wish to use the P\&R mode move from their origin to a P\&R facility and then transfer to public transportation to reach their final destination. In our model, all P\&R facilities share some common unobserved attributes other than location.

MNL and NL models are widely used to analyze transportation demand. Both models calculate the demand using the observed utility of alternatives ($V$), which generally consists of travel time and cost. \citet{holgui2012user} proposed a generalized utility for car and P\&R modes as follows:
\begin{align}
    \label{e:Vcar}
    V^\mathrm{c}_j & = c^\mathrm{time} TT^\mathrm{c}_j + c^\mathrm{cost} TC^\mathrm{c}_j + c^\mathrm{dist} TD^\mathrm{c}_j \\
    \label{e:Vpr}
    V^\mathrm{p}_{ij} & = c^\mathrm{time} TT^\mathrm{p}_{ij} + c^\mathrm{wait} WT^\mathrm{p}_{ij} + c^\mathrm{cost} TC^\mathrm{p}_{ij} + c^\mathrm{dist} TD^\mathrm{p}_{ij}.
\end{align}
$V^\mathrm{c}_j$ and $V^\mathrm{p}_{ij}$ denote the (observed) utility of using a private car and P\&R facility $i$, respectively, for O-D trip $j$; and $c^\mathrm{time}$, $c^\mathrm{cost}$, $c^\mathrm{dist}$, and $c^\mathrm{wait}$ are coefficients for the (in-vehicle) travel time, travel cost, travel distance, and wait ({\it i.e.}, out-of-vehicle) time, respectively. $TT^\mathrm{c}_j$, $TC^\mathrm{c}_j$, and $TD^\mathrm{c}_j$ denote the travel time, cost, and distance of using a private car for O-D trip $j$, and $TT^\mathrm{p}_{ij}$, $WT^\mathrm{p}_{ij}$, $TC^\mathrm{p}_{ij}$, and $TD^\mathrm{p}_{ij}$ denote the travel time, wait time, cost, and distance, respectively, of using P\&R facility $i$ for O-D trip $j$.

The P\&R FLP investigated in this study determines a fixed number ($N$) of optimal locations for P\&R facilities to be built. The value of $N$ is determined from the budget available for P\&R by dividing the available budget by the average cost of building a P\&R facility. 
The $N$ locations selected for the P\&R facilities will be indicated by binary variables $\left(x_i\in\{0,1\}:i\in\mathcal{P}\right)$ satisfying  
$$\sum_{i\in\mathcal{P}}x_i = N.$$
That is, $x_i = 1$ if facility $i\in \mathcal{P}$ is selected from among all of the candidates $\mathcal{P}$, and $x_i = 0$ otherwise.

The objective function of the P\&R FLP maximizes the demand for P\&R ({\it i.e.}, the total number of P\&R users). For an O-D trip $j\in\mathcal{T}$, the total number $R_j$ of travelers is given, and multiple P\&R facilities are available. The number of P\&R users for the O-D trip is the product of the total number of travelers on the O-D trip ($R_j$) and the probability of their choosing P\&R options for this O-D trip ($p_j^\mathrm{p}$). This study uses the NL model to calculate the probability of their choosing P\&R options, as the NL model considers the correlation between P\&R facilities. If there were no limits on the number of P\&R facilities to be built, the optimal combination of P\&R facilities would include all of the candidates $\mathcal{P}$. Because of the budget limitation, however, only a subset of the candidate P\&R facilities can eventually be built ({\it i.e.}, $N < |\mathcal{P}|$); thus, the demand for unselected P\&R candidate facilities must be excluded from the demand function for the expected total number of P\&R users. 

In this study, the NL demand function uses the nested structure shown in Fig.~\ref{f:nest}, where all P\&R options are in the same subtree. 
The first decision of a commuter for an O-D pair $j\in\mathcal{T}$ is whether to drive the car or take public transportation. It is represented by 
\begin{eqnarray}
    p^\mathrm{c}_j + p^\mathrm{p}_{j} = \frac{e^{V^\mathrm{c}_j}}{e^{V^\mathrm{c}_j}+e^{\lambda \Gamma^\mathrm{p}_j}} + \frac{e^{\lambda \Gamma^\mathrm{p}_j}}{e^{V^\mathrm{c}_j}+e^{\lambda \Gamma^\mathrm{p}_j}} = 1, ~\forall j \in \mathcal{T},\label{e:firstDecision}
\end{eqnarray}
where $\mathcal{T}$ denotes the set of all O-D trips, $\lambda$ is a logsum parameter, and $\Gamma^\mathrm{p}_j$ is the logsum of the exponential of the P\&R utilities for O-D trip $j$; {\it i.e.},
\begin{equation*}
    \Gamma^\mathrm{p}_j = \textrm{ln} \sum_{i \in \mathcal{P}} e^{V^\mathrm{p}_{ij}/\lambda} x_i, ~\forall j \in \mathcal{T}.
\end{equation*}

For O-D trip $j$, the two terms of (\ref{e:firstDecision}) are, respectively, the probability that the private car mode of transportation is chosen and the probability that the P\&R mode of transportation is chosen:  
\begin{eqnarray}
    \label{e:nlCar}
    p^\mathrm{c}_j &=& \frac{e^{V^\mathrm{c}_j}}{e^{V^\mathrm{c}_j}+e^{\lambda \Gamma^\mathrm{p}_j}}, ~\forall j \in \mathcal{T},\\
    \label{e:nlPR}
    p^\mathrm{p}_{j} &=& \frac{e^{\lambda \Gamma^\mathrm{p}_j}}{e^{V^\mathrm{c}_j}+e^{\lambda \Gamma^\mathrm{p}_j}}, ~\forall j \in \mathcal{T}.
\end{eqnarray}
If P\&R facility $i$ is not selected by the P\&R FLP, the utility of P\&R facility $i$ ({\it i.e.}, $V^\mathrm{p}_{ij}$) is not included for the calculation of $\Gamma^\mathrm{p}_j$. Specifically, the probability that a P\&R facility is chosen for O-D trip $j$ ({\it i.e.}, $p^\mathrm{p}_{j}$) excludes the utilities of the candidates $i\in\mathcal{P}$ for which $x_i=0$. 

Decomposing the probability ($p_j^{\mathrm{p}}$) that the P\&R option is chosen for O-D trip $j\in\mathcal{T}$, the decision-dependent probability that a specific P\&R facility $i\in\mathcal{P}$ is chosen is then represented by
\begin{equation}
    \label{e:nlTransit}
    p^\mathrm{p}_{ij} = p^\mathrm{p}_j \frac{e^{V^\mathrm{p}_{ij}/\lambda} x_i}{\sum_{k \in \mathcal{P}} e^{V^\mathrm{p}_{kj}/\lambda} x_k} = \frac{e^{\lambda \Gamma^\mathrm{p}_j}}{e^{V^\mathrm{c}_j}+e^{\lambda \Gamma^\mathrm{p}_j}} \frac{e^{V^\mathrm{p}_{ij}/\lambda} x_i}{\sum_{k \in \mathcal{P}} e^{V^\mathrm{p}_{kj}/\lambda} x_k}, ~\forall i \in \mathcal{P}, ~\forall j \in \mathcal{T}.
\end{equation}
Thus, $p^\mathrm{p}_{ij}$ is nonzero if and only if $x_i$ is not zero, and the objective function of the P\&R FLP is the demand function
$$\sum_{i \in \mathcal{P}} \sum_{j \in \mathcal{T}} R_j p^\mathrm{p}_{ij}
=\sum_{j \in \mathcal{T}} R_j p^\mathrm{p}_{j},$$ 
which is the expected total number of P\&R users.

Now, the P\&R FLP under the NL demand function is formalized as a nonlinear mathematical programming formulation as follows:
\begin{align}
    \label{e:obj}
    \textrm{Maximize } & \sum_{i \in \mathcal{P}} \sum_{j \in \mathcal{T}} R_j p^\mathrm{p}_{ij} \\
    \textrm{subject to} & \nonumber \\
    \label{e:selection}
    & \sum_{i \in \mathcal{P}} x_i = N\\
    \label{e:nl1}
    & \Gamma^\mathrm{p}_j = \textrm{ln} \sum_{i \in \mathcal{P}} e^{V^\mathrm{p}_{ij}/\lambda} x_i, ~\forall j \in \mathcal{T} \\
    \label{e:nl2}
    & p^\mathrm{p}_{ij} = \frac{e^{\lambda \Gamma^\mathrm{p}_j}}{e^{V^\mathrm{c}_j}+e^{\lambda \Gamma^\mathrm{p}_j}} \frac{e^{V^\mathrm{p}_{ij}/\lambda} x_i}{\sum_{k \in \mathcal{P}} e^{V^\mathrm{p}_{kj}/\lambda} x_k}, ~\forall i \in \mathcal{P}, ~\forall j \in \mathcal{T} \\
    \label{e:var}
    & x_i \in \{ 0, 1 \}, ~\forall i \in \mathcal{P}.
\end{align}
Equation~(\ref{e:selection}) ensures that only $N$ park-and-ride facilities are to be built.
{As (\ref{e:nl1}) and (\ref{e:nl2}) are functions 
of $x=\left(x_i:i\in\mathcal{P}\right)$ and can be substituted for $p^{\mathrm{p}}_{ij}=p^{\mathrm{p}}_{ij}(x)$ in the objective function, Equation~(\ref{e:selection}) is the only constraint of the formulation.
If $\lambda\neq 1$, Equations~(\ref{e:nl1}) and (\ref{e:nl2}) do not allow linearization of the nonlinear program.} 

\subsection{Relationship to Multinomial Logit Model}
The MNL model is a special case of our NL model. The logsum parameter $\lambda$ of the NL model is bounded by zero and one. When $\lambda=1$, P\&R facilities are uncorrelated ({\it i.e.}, there is no nested structure), and 
$$e^{\lambda \Gamma^\mathrm{p}_j}=\sum_{i \in \mathcal{P}} e^{V^\mathrm{p}_{ij}} x_i.$$ 
Then, (\ref{e:nlCar}) and (\ref{e:nlTransit}) are
\begin{align}
    \label{e:mnlCar}
    p^\mathrm{c}_j & = \frac{e^{V^\mathrm{c}_j}}{e^{V^\mathrm{c}_j}+\sum_{i \in \mathcal{P}} e^{V^\mathrm{p}_{ij}} x_i}, ~\forall j \in \mathcal{T}\\
    \label{e:mnlTransit}
    p^\mathrm{p}_{ij} &= \frac{e^{V^\mathrm{p}_{ij}} x_i}{e^{V^\mathrm{c}_j}+\sum_{k \in \mathcal{P}} e^{V^\mathrm{p}_{kj}} x_k}, ~\forall i \in \mathcal{P}, ~\forall j \in \mathcal{T}.
\end{align}

Then, the IIA condition between private car and P\&R facility $i$ for O-D trip $j$ is 
\begin{equation}
    \label{e:iia}
    \frac{p^\mathrm{p}_{ij}}{p^\mathrm{c}_j} = \frac{e^{V^\mathrm{p}_{ij}} x_i}{e^{V^\mathrm{c}_j}}, ~\forall i \in \mathcal{P}, ~\forall j \in \mathcal{T},
\end{equation}
{by dividing (\ref{e:mnlTransit}) by (\ref{e:mnlCar}). 
The IIA condition means that the ratio of $p^\mathrm{p}_{ij}$ to $p^\mathrm{c}_j$ depends only on the utilities of P\&R facility $i$ with $x_i=1$ ({\it i.e.}, a facility to be built) and of the private car.} Note that (\ref{e:mnlCar}) and (\ref{e:mnlTransit}) are equivalent to (\ref{e:iia}) taken together with 
\begin{eqnarray*}
p^\mathrm{c}_j + \sum_{i \in \mathcal{P}} p^\mathrm{p}_{ij} = 1.
\end{eqnarray*}
Similarly, the IIA condition between P\&R facilities $i$ and $k$ for O-D trip $j$ is
\begin{equation}
    \label{e:iia2}
    \frac{p^\mathrm{p}_{kj}}{p^\mathrm{p}_{ij}} = \frac{e^{V^\mathrm{p}_{kj}} x_k}{e^{V^\mathrm{p}_{ij}} x_i}, ~\forall i, k \in \mathcal{P}, ~\forall j \in \mathcal{T}.
\end{equation}

Unlike the P\&R FLP with the NL model presented in (\ref{e:obj})--(\ref{e:var}), the P\&R FLP with the MNL model can be expressed as an MILP formulation. 
Of the three linearized formulations presented by \citet{haase2014comparison}, the formulation of \citet{aros2013p} is adopted for this study because it is quickly solved. The MILP formulation of the P\&R FLP with the MNL model is 
\begin{align}
    \label{e:obj2}
    \textrm{Maximize } & \sum_{i \in \mathcal{P}} \sum_{j \in \mathcal{T}} R_j p^\mathrm{p}_{ij} \\
    \label{e:selection2}
    \textrm{subject to} & \nonumber \\
    & \sum_{i \in \mathcal{P}} x_i = N\\
    \label{e:probConst}
    & p^\mathrm{p}_{ij} \leq x_i, ~\forall i \in \mathcal{P}, ~\forall j \in \mathcal{T} \\
    \label{e:totalConst}
    & p^\mathrm{c}_j + \sum_{i \in \mathcal{P}} p^\mathrm{p}_{ij} = 1, ~\forall j \in \mathcal{T}\\
    \label{e:iiaConst2}
    & p^\mathrm{p}_{ij} \leq \frac{e^{V^\mathrm{p}_{ij}}}{e^{V^\mathrm{c}_j}} p^\mathrm{c}_j, ~\forall i \in \mathcal{P}, ~ \forall j \in \mathcal{T}  \\
    \label{e:iiaConst3}
    & p^\mathrm{c}_j \leq \frac{e^{V^\mathrm{c}_j}}{e^{V^\mathrm{p}_{ij}}} p^\mathrm{p}_{ij} + (1-x_i), ~\forall i \in \mathcal{P}, ~\forall j \in \mathcal{T} \\
    \label{e:iiaConst1}
    & p^\mathrm{p}_{ij} \leq \frac{e^{V^\mathrm{p}_{ij}}}{e^{V^\mathrm{p}_{kj}}} p^\mathrm{p}_{kj} + (1-x_k), ~\forall i, k \in \mathcal{P}, ~ \forall j \in \mathcal{T} \\
    & x_i \in \{ 0, 1 \}, ~\forall i \in \mathcal{P} \\
    & p^\mathrm{c}_j \geq 0, ~\forall j \in \mathcal{T} \\
    \label{e:var2}
    & p^\mathrm{p}_{ij} \geq 0, ~\forall i \in \mathcal{P}, ~ \forall j \in \mathcal{T}.
\end{align}

Equations~(\ref{e:mnlCar}) and (\ref{e:mnlTransit}) are formalized by the linear constraints~(\ref{e:probConst})--(\ref{e:iiaConst3}).
Equation~(\ref{e:probConst}) forces $p^\mathrm{p}_{ij} = 0$ if $x_i = 0$, as given in (\ref{e:mnlTransit}). 
Equation~(\ref{e:totalConst}) ensures that the sum of all probabilities ({\it e.g.}, $p^\mathrm{c}_j$ and $p^\mathrm{p}_{ij}$) is equal to one for each O-D trip. 
The IIA condition of (\ref{e:iia}) with $x_i=1$ is the equation induced by (\ref{e:iiaConst2}) and (\ref{e:iiaConst3}).
Similarly, the IIA condition of (\ref{e:iia2}) is ensured by (\ref{e:iiaConst1}).
{
Note that Equation~(\ref{e:iiaConst1}) is redundant as (\ref{e:mnlCar}) and (\ref{e:mnlTransit}) are implied by (\ref{e:probConst})--(\ref{e:iiaConst3}).
In a similar manner, we may alternatively delete Equations~(\ref{e:iiaConst2}) and (\ref{e:iiaConst3}) instead of Equation~(\ref{e:iiaConst1}) (see \citet{haase2014comparison}).
To further strengthen the formulation, \citet{krohn2021} have introduced a tighter bound on $p^{\mathrm{p}}_{ij}$ in (\ref{e:probConst}).
}

\section{{Neighborhood Search Method}}\label{s:neighbor}
In this section, we develop a neighborhood search (NS) method to solve the P\&R FLP by considering the polyhedron whose vertices are the feasible solutions $x^{\mathrm{INT}}\in\{0,1\}^{|\mathcal{P}|}$ satisfying (\ref{e:selection}). 
We do not consider the entire polyhedron but rather the adjacency graph defined by the vertices and edges of the polyhedron. 
Like the simplex method, the neighborhood search method starts at one vertex of the polyhedron and moves along an edge to the next vertex (called a neighbor), improving the demand function. 
The method iterates until there is no better neighbor at a local maximum.

Let $\mathrm{INT}(N)$ denote the set of $x^{\mathrm{INT}}\in\{0,1\}^{|\mathcal{P}|}$ satisfying (\ref{e:selection}), which are the feasible solutions to an $N$-hub facility location problem.
Let $\mathrm{DEMAND}_\lambda(x)$ denote the objective function (\ref{e:obj}) in which $p_{ij}^{\mathrm{p}}=p_{ij}^{\mathrm{p}} (x)$ are functions of $x$ defined by (\ref{e:nl1}) and (\ref{e:nl2}).
More precisely, $p_{ij}^{\mathrm{p}}=p_{ij}^{\mathrm{p}} (x)$ is determined by variable vector $x$ along with parameter $\lambda$, and therefore $\mathrm{DEMAND}_\lambda(x)$ is determined by $x$ and $\lambda$. We will denote $\mathrm{DEMAND}_\lambda(x)$ simply as $\mathrm{DEMAND}(x)$ when $\lambda$ does not need to be specified. 
Then, the P\&R FLP can be simply written as 
\begin{eqnarray}
\max\left\{\mathrm{DEMAND}\left(x^{\mathrm{INT}}\right):x^{\mathrm{INT}}\in\mathrm{INT}(N)\right\}.\label{e:brief}
\end{eqnarray}
The relaxation of $\mathrm{INT}(N)$ is 
$$\mathrm{REL}(N) = \left\{ \tilde{x}\in [0,1]^{|\mathcal{P}|} : \tilde{x}\mbox{ satisfying (\ref{e:selection})} \right\},$$
where $[0,1]$ is the interval $\{\tilde{x}_i\in\mathbb{R}:0\leq \tilde{x}_i\leq 1\}$ between 0 and 1.
It is easy to see that $\mathrm{INT}(N)$ is the set consisting of the vertices of the polyhedron $\mathrm{REL}(N)$:
\begin{proposition}
The feasible solutions $x^{\mathrm{INT}}\in\mathrm{INT}(N)$ are the vertices of the polyhedron $\mathrm{REL}(N)$.
\end{proposition}
Considering $\mathrm{INT}(N)$ as a subset of $\mathrm{REL}(N)$, we will refer to the feasible solutions $x^{\mathrm{INT}}$ interchangeably as the integer solutions, particularly when we distinguish the integer solutions from the fractional solutions $\tilde{x}\in\mathrm{REL}(N)\setminus \mathrm{INT}(N)$. 

It is also easy to see that two vertices $x^{\mathrm{INT}(1)}$ and $x^{\mathrm{INT}(2)}$ are adjacent ({\it i.e.}, connected by an edge of the polyhedron $\mathrm{REL}(N)$) if their 1-norm distance is 2: 
\begin{proposition}
Two vertices $x^{\mathrm{INT}(1)}$ and $x^{\mathrm{INT}(2)}$ are adjacent on the polyhedron $\mathrm{REL}(N)$ if and only if
$$\left\| x^{\mathrm{INT}(1)} - x^{\mathrm{INT}(2)} \right\|_1 = \sum_{i\in\mathcal{P}} \left|x^{\mathrm{INT}(1)}_i - x^{\mathrm{INT}(2)}_i\right| = 2.$$
\end{proposition}
When two vertices $x^{\mathrm{INT}(1)}$ and $x^{\mathrm{INT}(2)}$ are adjacent, one is said to be a \emph{neighbor} of the other.
Note that each vertex has $N\times (|\mathcal{P}|-N)$ neighbors.
On the adjacency graph $G=(V, E)$ of the polyhedron $\mathrm{REL}(N)$, where $V$ is the set of the vertices and $E$ is the set of the edges of the polyhedron, a \emph{local optimum} $x^{\mathrm{LOPT}}$ is defined to have all inferior neighbors $x^{\mathrm{NBR}(k)}$; {\it i.e.}, $$\mathrm{DEMAND}\left( x^{\mathrm{LOPT}}\right)\geq\mathrm{DEMAND}\left( x^{\mathrm{NBR}(k)}\right) \mbox{ for } k=1,...,N\times (|\mathcal{P}|-N).$$
For further details regarding the vertices and edges of a polyhedron, the reader may refer to \citet{bertsimas1997}.

\begin{figure}
    \centering
    \includegraphics[width=0.7\textwidth]{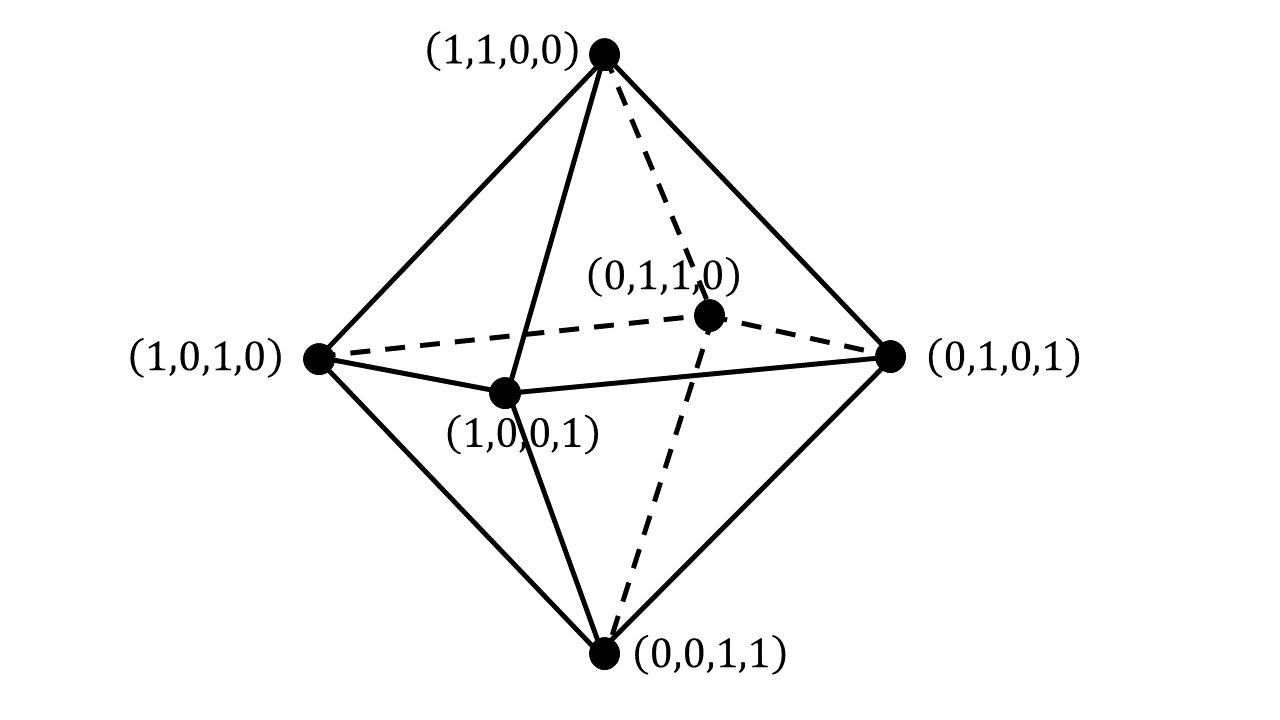}
    \caption{Adjacency graph of $\mathrm{REL}(N)$ selecting $N=2$ locations among $|\mathcal{P}|=4$ candidates}
    \label{f:octa}
\end{figure}

{
Since $\dot{x}=\left(\dot{x}_i=N/|\mathcal{P}|\mbox{ for }i\in \mathcal{P}\right)$ is an inner point that satisfies Equation~(\ref{e:selection}) and the strict inequalities $0< \dot{x}_i = N/|\mathcal{P}|<1$ for $i\in \mathcal{P}$, the dimension of the polyhedron $\mathrm{REL}(N)$ is $|\mathcal{P}|-1$. In particular, $\dot{x}$ is the average of the integer solutions $x^{\mathrm{INT}}\in\mathrm{INT}(N)$:
\begin{proposition} 
The average of the integer solutions $x^{\mathrm{INT}}\in\mathrm{INT}(N)$ is $\dot{x}$; {\it i.e.},
\begin{equation*}
\dot{x} = \frac{\sum\left\{x^{\mathrm{INT}}\in\mathrm{INT}(N)\right\} }{\left|\mathrm{INT}(N)\right|}.
\end{equation*}
\end{proposition}
\begin{pf}
Each facility $i\in\mathcal{P}$ is a member of the same number $\binom{|\mathcal{P}|-1}{N-1}$ of combinations of $N$ facilities. Then,
\begin{eqnarray*}
\dot{x}_i = \frac{N}{|\mathcal{P}|} = \frac{\binom{|\mathcal{P}|-1}{N-1}}{\binom{|\mathcal{P}|}{N}} = \frac{\binom{|\mathcal{P}|-1}{N-1}}{\left|\mathrm{INT}(N)\right|} = \frac{\sum_{x^{\mathrm{INT}}\in\mathrm{INT}(N)} x^{\mathrm{INT}}_i}{\left|\mathrm{INT}(N)\right|}.
\end{eqnarray*}
\qed
\end{pf}
We refer to $\dot{x}$ as the \emph{centroid} of $\mathrm{REL}(N)$ as it is the center of weight of the polyhedron $\mathrm{REL}(N)$. For further details regarding inner points and the dimensionality of a polyhedron, the reader may refer to \citet{nemhauser1999}. 
}

Fig.~\ref{f:octa} depicts the adjacency graph $G=(V,E)$ of $\mathrm{REL}(N)$ selecting $N=2$ locations from among $|\mathcal{P}|=4$ candidates. The graph has $|V|=6$ vertices and $|E|=12$ edges. Each vertex is adjacent to $N\times (|\mathcal{P}|-N)=4$ neighbors. Two vertices are adjacent whenever exactly two components differ. The centroid is $\dot{x}=(0.5,0.5,0.5,0.5)$. The polyhedron $\mathrm{REL}(N)$ is three-dimensional as $|\mathcal{P}|-1=3$. 

\begin{algorithm}
\caption{Neighborhood Search Method}
\label{a:deterministic}
\begin{algorithmic}[1]
\STATE Start with any integer solution $x^{\mathrm{INT}}$\label{l:start}
\STATE Set $\sf is\_local\_optimum = False$
\WHILE{$\sf is\_local\_optimum = False$}\label{l:whileStarts}
    \STATE Set $\sf is\_local\_optimum = True$
    \FOR{$i\in\mathrm{supp}\left(x^{\mathrm{INT}}\right)$ and $k\in\mathcal{P}\setminus\mathrm{supp}\left(x^{\mathrm{INT}}\right)$}\label{l:forStarts}
    \STATE Set $x^{\mathrm{NBR}}\leftarrow x^{\mathrm{INT}}$
    \STATE $x^{\mathrm{NBR}}_i\leftarrow 0$
    \STATE $x^{\mathrm{NBR}}_k\leftarrow 1$
    \IF{$\mathrm{DEMAND}\left(x^{\mathrm{INT}}\right) < \mathrm{DEMAND}\left(x^{\mathrm{NBR}}\right)$}\label{l:betterNeighbor}
    \STATE Update $x^{\mathrm{INT}}\leftarrow x^{\mathrm{NBR}}$ \label{l:update}
    \STATE Set $\sf is\_local\_optimum = False$ \label{l:move}
    \STATE {\bf break}
    \ENDIF
    \ENDFOR \label{l:forEnds}
\ENDWHILE \label{l:whileEnds}
\RETURN the local maximum $x^{\mathrm{LOPT}}=x^{\mathrm{INT}}$
\end{algorithmic}
\end{algorithm}

Pseudocode for the neighborhood search method is given in Algorithm~\ref{a:deterministic}. In Line~\ref{l:start}, the search starts with an integer solution $x^{\mathrm{INT}}$, which is the initial best known integer solution.
In Lines~\ref{l:whileStarts}--\ref{l:whileEnds}, the method moves the best known integer solution $x^{\mathrm{INT}}$ to its neighbor, improving $\mathrm{DEMAND}\left(x^{\mathrm{INT}}\right)$, until $x^{\mathrm{INT}}$ reaches a local optimum $x^{\mathrm{LOPT}}$. During this process, $\sf is\_local\_optimum = False$ indicates that $x^{\mathrm{INT}}$ is not guaranteed to be a local optimum, in which case it continues moving.
In Lines~\ref{l:forStarts}-\ref{l:forEnds}, the method enumerates the neighbors $x^{\mathrm{NBR}}$ of $x^{\mathrm{INT}}$. If a superior neighbor is found (Line~\ref{l:betterNeighbor}),
the method updates the best known integer solution with the superior neighbor by setting $x^{\mathrm{INT}}$ to $x^{\mathrm{NBR}}$ (Line~\ref{l:update}) and sets $\sf is\_local\_optimum$ to $\sf False$ so that the best known integer solution will continue moving (Line~\ref{l:move}).
If there is no superior neighbor ({\it i.e.}, $\sf is\_local\_optimum=True$), the method returns the local optimum $x^{\mathrm{LOPT}}=x^{\mathrm{INT}}$.
Note that the 1-norm distance between $x^{\mathrm{INT}}$ and $x^{\mathrm{NBR}}$ equals 2; {\it i.e.}, $x^{\mathrm{NBR}}$ is given by swapping a location $i\in\mathrm{supp}\left( x^{\mathrm{INT}}\right)=\left\{ i\in\mathcal{P}:x^{\mathrm{INT}}_i=1\right\}$ of the support of $x^{\mathrm{INT}}$ and a non-support location $k\in\mathcal{P}\setminus\mathrm{supp}\left(x^{\mathrm{INT}}\right)$. 


As will be shown in Section~\ref{sec:experiments}, computational experiments verified that the neighborhood search method solves the medium-scale P\&R FLP to exact optimality extremely rapidly.
Specifically, it solves the MNL model of the medium-scale instances to exact optimality 10,000 times faster than the MILP formulation introduced by \citet{aros2013p}. 
On the large-scale instances, which the MILP formulation cannot solve, the neighborhood search method solves the MNL model within 30 s.

\section{Adaptive Randomized Rounding Procedure}\label{sec:heuristics}
In this section, the neighborhood search method described in the previous section is randomized and developed into an adaptive randomized rounding (ARR) procedure.
Given a fractional solution $\tilde{x}$, which we will refer to as a \emph{seed}, randomized rounding finds an integer solution $x^{\mathrm{INT}}$ near the seed. A randomized rounding procedure finds multiple integer solutions and selects the best known integer solution $x^{\mathrm{BKI}}$. An adaptive randomized rounding procedure moves the seed $\tilde{x}$ toward the best known integer solution $x^{\mathrm{BKI}}$ to find a better integer solution near $x^{\mathrm{BKI}}$ ({\it i.e.}, a better integer solution in the approximate neighborhood of $x^{\mathrm{BKI}}$).

The traditional randomized rounding procedure first solves the relaxation 
\begin{eqnarray}
\mathrm{DEMAND} (\tilde{x}^*) = \max\left\{ \mathrm{DEMAND} (\tilde{x}) : \tilde{x}\in \mathrm{REL} (N)\right\}\label{e:relax}
\end{eqnarray}
and then finds an integer solution $X^{\mathrm{INT}}=x^{\mathrm{INT}}$ based on the fractional optimal solution $\tilde{x}^*$ such that
\begin{eqnarray}
P\left[X^{\mathrm{INT}}_i=1\right]=\tilde{x}^*_i\mbox{ for }i\in \mathcal{P},\label{e:probability}
\end{eqnarray}
where $X$ denotes the random variable for an integer solution $x$.
{
(For further details regarding (traditional) randomized rounding, the readers may refer to \citet{motwani1995} and \citet{vazirani2001}.)
}
A randomized rounding procedure conducts multiple trials to find integer solutions $x^{\mathrm{INT}(t)},t=1,...,T$, and selects the solution found (BFS); {\it i.e.}, 
$$\mathrm{BSF}=\arg\max\left\{ \mathrm{DEMAND} \left( x^{\mathrm{INT}(t)}\right): t=1,...,T\right\}.$$

As it is not easy to solve the relaxation (\ref{e:relax}), the initial seed for our randomized rounding procedure will be the uniform distribution $\dot{x}=\left( \dot{x}_i = N/|\mathcal{P}| : i\in\mathcal{P}\right)$, and the seed will move to converge to the best known integer solution $x^{\mathrm{BKI}}$ so far, searching for better integer solutions in the vicinity of the best known integer solution. 
Because $\tilde{x}$ is interpreted as the likelihood of $x^{\mathrm{INT}}$, we consider only $\tilde{x}\in [0,1]^{|\mathcal{P}|}\subseteq\mathbb{R}^{|\mathcal{P}|}$, with the components $\tilde{x}_i$ being between 0 and 1, throughout this study.
Section~\ref{s:three} describes the three main components of our adaptive randomized rounding technique, and Section~\ref{s:ARR} shows the construction of the full adaptive randomized rounding procedure.


\subsection{Three Main Components}\label{s:three}

\subsubsection{Rounding}\label{s:rounding}
\begin{figure}
    \centering
    \includegraphics[width=0.7\textwidth]{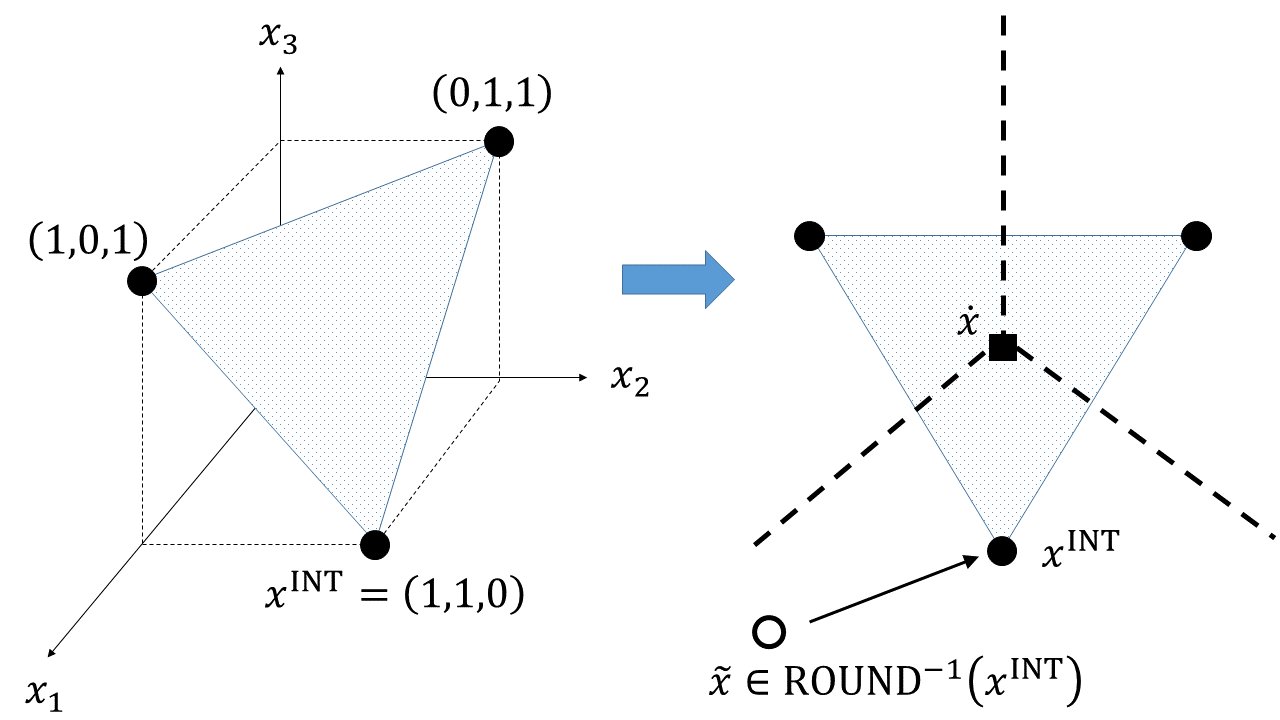}
    \caption{The hyperplane $H=\left\{\tilde{x}\in\mathrm{R}^3:\tilde{x}_1+\tilde{x}_2+\tilde{x}_3=2\right\}$ is partitioned into three regions $H\cap \mathrm{ROUND}^{-1}\left(x^{\mathrm{INT}}\right)$ corresponding to the three integer solutions $x^{\mathrm{INT}}\in\mathrm{INT}(N)$ where $N=2$ and $|\mathcal{P}|=3$. The rounding function maps a fractional solution $\tilde{x}$ in the region of $\mathrm{ROUND}^{-1}\left(x^{\mathrm{INT}}\right)\subseteq [0,1]^{|\mathcal{P}|}$ to $x^{\mathrm{INT}}$.}
    \label{f:regions}
\end{figure}

A rounding function $\mathrm{ROUND}:\mathbb{R}^{|\mathcal{P}|}\rightarrow\mathrm{INT}(N)$ maps a fractional solution $\tilde{x}\in\mathbb{R}^{|\mathcal{P}|}$ to the nearest integer solution $x^{\mathrm{INT}}=\mathrm{ROUND}\left(\tilde{x}\right)$ indicating the $N$ largest components; {\it i.e.}, $x^{\mathrm{INT}}_i=1$, if and only if $\tilde{x}_i$ is one of the $N$ largest components of $\tilde{x}$. (A tie may be broken by a linear ordering $\prec$ of $\mathcal{P}$. That is, the same $N$-th largest components $\tilde{x}_i=\tilde{x}_k$ of $i\prec k\in\mathcal{P}$ may be rounded to $x^{\mathrm{INT}}_i=0$ and $x^{\mathrm{INT}}_k=1$.) Then, an integer solution $x^{\mathrm{INT}}$ is rounded to itself; {\it i.e.}, 
$$\mathrm{ROUND}\left(x^{\mathrm{INT}}\right)=x^{\mathrm{INT}}.$$ 
The inverse image of the rounding function partitions the entire space $\mathbb{R}^{|\mathcal{P}|}$ into regions $\mathrm{ROUND}^{-1} \left( x^{\mathrm{INT}}\right)$ for $x^{\mathrm{INT}}\in\mathrm{INT}(N)$.
In Fig.~\ref{f:regions}, the regions $\left\{\mathrm{ROUND}^{-1} \left( x^{\mathrm{INT}}\right):x^{\mathrm{INT}}\in\mathrm{INT}(N)\right\}$ partition the hyperplane $H$ defined by (\ref{e:selection}), where $N=2$ and $\mathcal{P}=\{ 1,2,3\}$. The rounding function maps a fractional solution $\tilde{x}$ lying in the region $\mathrm{ROUND}^{-1}\left(x^{\mathrm{INT}}\right)$ to $x^{\mathrm{INT}}$. 

\subsubsection{Randomization}
Our randomized rounding procedure perturbs a given seed $\tilde{x}$ randomly so that the perturbed seed $\tilde{x}^{\mathrm{PTB}}$ may lie in $\mathrm{ROUND}^{-1}(x^{\mathrm{INT}})\subseteq [0,1]^{|\mathcal{P}|}$, satisfying (\ref{e:proportional}).
A random perturbation increases each component $\tilde{x}_i$ by uniform distribution $\mathrm{U}\left[ \tilde{x},1\right]$ in the interval $\left[ \tilde{x},1\right]$; {\it i.e.},
\begin{eqnarray}
\tilde{x}^{\mathrm{PTB}}_i = \tilde{x}_i + \left( 1 - \tilde{x}_i\right)\tilde{\pi}_i\mbox{ with }\tilde{\pi}_i\in \mathrm{U}\left[ 0,1\right]\mbox{ for }i\in\mathcal{P}.\label{e:perturb}
\end{eqnarray}
{
If $\tilde{x}_i$ equals 1, $\tilde{x}^{\mathrm{PTB}}_i (\tilde{x}) = 1$ is always the highest component, and facility $i$ is selected ({\it i.e.}, $x^{\mathrm{INT}}_i=1$). 
Therefore, if the seed is an integer solution $\tilde{x}=x^{\mathrm{INT}}$, the perturbed solution lies in $\tilde{x}^{\mathrm{PTB}}\in\mathrm{ROUND}^{-1}\left( x^{\mathrm{INT}}\right)$; {\it i.e.},
}
{
\begin{eqnarray}
\mathrm{ROUND}\left( \tilde{x}^{\mathrm{PTB}}\right)=x^{\mathrm{INT}}\mbox{ whenever }\tilde{x}=x^{\mathrm{INT}}.
\end{eqnarray}
}
In contrast, $x^{\mathrm{INT}}_i$ does not need to be zero for a facility candidate $i$ with $\tilde{x}_i = 0$, giving the facility candidate a positive opportunity ({\it i.e.}, $\tilde{x}^{\mathrm{PTB}}>0$) to be selected, whenever $\tilde{x}$ is not an integer solution. In the remainder of this paper, we denote $x^{\mathrm{INT}} = \mathrm{ROUND} \left( \tilde{x}^{\mathrm{PTB}} \right)$ by $x^{\mathrm{INT}}\left( \tilde{x} \right)$.

\subsubsection{Moving Seed}\label{s:moving}
As it is not easy to solve the relaxation (\ref{e:relax}), the initial seed for our randomized rounding procedure will be the uniform distribution $\tilde{x}^{(0)}=\dot{x}=\left( \dot{x}_i = N/|\mathcal{P}| : i\in\mathcal{P}\right)$, and the seed will move to converge to the best known integer solution $x^{\mathrm{BKI}}$ as better integer solutions are sought in the vicinity of $x^{\mathrm{BKI}}$. The procedure first finds an integer solution $x^{\mathrm{INT}(0)}=\mathrm{ROUND}\left( \tilde{x}^{\mathrm{PTB}(0)}\right)$, where $\tilde{x}^{\mathrm{PTB}(0)}$ is given by perturbing $\tilde{x}^{(0)}=\dot{x}$. The first integer solution is the best known integer solution $x^{\mathrm{BKI}}=x^{\mathrm{INT}(0)}$. At each trial $t\geq 1$, the procedure then moves the last seed $\tilde{x}^{(t-1)}$, converging to the best known integer solution $x^{\mathrm{BKI}}$ by exponential smoothing:
\begin{equation}
\tilde{x}^{(t)} = \left( 1 - \alpha (t) \right)\tilde{x}^{(t-1)} + \alpha (t) x^{\mathrm{BKI}},\label{eqn:expSmoothNew}
\end{equation}
where $0<\alpha (t) < 1$. 
When $\mathrm{DEMAND}\left( x^{\mathrm{BKI}}\right) < \mathrm{DEMAND}\left( x^{\mathrm{INT}(t)}\right)$, the procedure updates the best known integer solution $x^{\mathrm{BKI}}$ to $x^{\mathrm{INT}(t)}$.

{
The initial seed $\tilde{x}^{(0)}=\dot{x}$ is the centroid, which is the probability distribution of $X^{\mathrm{INT}}\left(\dot{x}\right)$:
\begin{proposition} 
At the centroid $\dot{x}$, the probability distribution is $\dot{x}$ itself; {\it i.e.},
$$P\left[X_i^{\mathrm{INT}}\left(\dot{x}\right)=1\right] = \dot{x}_i.$$
\end{proposition}
}
\begin{pf}
{
On the seed of $\dot{x}$, each candidate $i\in\mathcal{P}$ can be selected equally likely. Therefore,
$$P\left[X_i^{\mathrm{INT}}\left(\dot{x}\right)=1\right] = \frac{\sum_{i\in\mathcal{P}}P\left[X_i^{\mathrm{INT}}\left(\dot{x}\right)=1\right]}{|\mathcal{P}|}
= \frac{\sum_{i\in\mathcal{P}}E\left[X_i^{\mathrm{INT}}\left(\dot{x}\right)\right]}{|\mathcal{P}|}
= \frac{E\left[\sum_{i\in\mathcal{P}} X_i^{\mathrm{INT}}\left(\dot{x}\right)\right]}{|\mathcal{P}|}=\frac{N}{|\mathcal{P}|}=\dot{x}_i.$$
\qed
}
\end{pf}
{
That is, the centroid $\dot{x}$ satisfies Equation~(\ref{e:probability}) in replace of the optimal solution $\tilde{x}^*$ to the relaxation~(\ref{e:relax}). At the centroid, the randomized rounding procedure finds all integer solutions equally likely. 
}

{
Beginning with the centroid $\tilde{x}^{(0)}=\dot{x}$, the seed $\tilde{x}^{(t)}$ (at trial $t\geq 1$) will move to converge to the best known integer solution by exponential smoothing~(\ref{eqn:expSmoothNew}), searching for better integer solutions in the vicinity of the best known integer solution $x^{\mathrm{BKI}}$. Instead of strictly satisfying (\ref{e:probability}), the probability distribution will be proportional to the seed $\tilde{x}$; more precisely, 
\begin{eqnarray}
\tilde{x}_i<\tilde{x}_k\Rightarrow P\left[ X^{\mathrm{INT}}_i = 1\right] < P\left[ X^{\mathrm{INT}}_k = 1\right]\mbox{ for }i\neq k\in\mathcal{P}.\label{e:proportional}
\end{eqnarray}
Converging toward $x^{\mathrm{BKI}}$, the seed $\tilde{x}^{(t)}$ finds the same integer solution $x^{\mathrm{INT}(t)}=x^{\mathrm{INT}}\left(\tilde{x}^{(t)}\right)=x^{\mathrm{BKI}}$ increasingly frequently. When the seed is too close to the best known integer solution, our randomized rounding procedure cannot find a better integer solution and  will reset the seed back to the centroid $\dot{x}$.
}

\subsection{Full Procedure}\label{s:ARR}
In this section, the full adaptive randomized rounding (ARR) procedure is developed. For this, we will define the rule for slowing the convergence speed ({\it i.e.}, smoothing constant) $\alpha (t)$ in (\ref{eqn:expSmoothNew}) and the rule for resetting the seed.

The initial seed for our adaptive randomized rounding procedure is the centroid $\tilde{x}^{(0)}=\dot{x}$. The centroid finds the initial integer solution $x^{\mathrm{BKI}}=x^{\mathrm{INT}}\left(\dot{x}\right)$. Each trial $t\geq 1$ moves the seed, converging to the best known integer solution $x^{\mathrm{BKI}}$ by exponential smoothing (\ref{eqn:expSmoothNew}). When $\tilde{x}^{(t)}$ is too close to $x^{\mathrm{BKI}}$, a run of many trials will find the same integer solution $x^{\mathrm{INT}(t)}=x^{\mathrm{BKI}}$, in which case we reset the seed $\tilde{x}^{(t)}$ back to the centroid $\dot{x}$. 

For computational convenience, we further relax the relaxation $\mathrm{REL}(N)$ to $\left[ 0,1\right]^{|\mathcal{P}|}$ and consider its centroid $\ddot{x}$ given by
\begin{eqnarray}
\ddot{x}_i=0.5\mbox{ for }i\in\mathcal{P}.
\end{eqnarray}
This is equivalent to $\dot{x}=\left(\dot{x}_i=N/|\mathcal{P}|\right)$ in that 
$$P\left[X^{\mathrm{INT}}_i \left( \ddot{x}\right)=1\right] = P\left[X^{\mathrm{INT}}_i \left( \dot{x}\right)=1\right] = N/|\mathcal{P}|.$$
In general, $\tilde{x}'\in\left[ 0,1\right]^{|\mathcal{P}|}$ is projected to $\tilde{x}\in\left[ 0,1\right]^{|\mathcal{P}|}$, satisfying $\sum_{i\in\mathcal{P}}\tilde{x}_i=N$ as follows.
\begin{theorem}
Let $\tilde{x}'\in [0,1]^{|\mathcal{P}|}$, and let
\begin{equation}
    \tilde{x} = \mathbf{1}- \frac{|\mathcal{P}| - N}{\lVert \mathbf{1} - \tilde{x}'\rVert_1}\left( \mathbf{1} - \tilde{x}' \right).
\end{equation}
Then  
\begin{equation}
    P\left[ X^{\mathrm{INT}} \left( \tilde{x}'\right) = x^{\mathrm{INT}}\right]
    =     P\left[ X^{\mathrm{INT}} \left( \tilde{x}\right) = x^{\mathrm{INT}}\right]
    \mbox{ for }x^{\mathrm{INT}}\in\mathrm{INT}(N).
\end{equation}
\end{theorem}

\begin{pf}
An integer solution $x^{\mathrm{INT}}\left(\tilde{x}'\right)$ selects the largest $N$ values $\tilde{l}'_i$ of
$$\tilde{l}'_i = \tilde{x}'_i + \left( 1 - \tilde{x}'_i \right) \tilde{\pi}_i$$
with $\tilde{\pi}_i\in \mathrm{U} \left[ 0,1 \right]$ for $i\in\mathcal{P}$.
This is equivalent to selecting the smallest $N$ values of
\begin{equation}
    1 - \tilde{l}'_i = \left(1-\tilde{x}'_i\right) \left( 1 - \tilde{\pi}_i\right) = \left( \frac{\lVert \mathbf{1}-\tilde{x}'\rVert_1}{ |\mathcal{P}| - N } \right) \left(1-\tilde{x}_i\right) \left( 1 - \tilde{\pi}_i\right)\mbox{ for }i\in\mathcal{P}.\label{e:transform}
\end{equation}
Since $\left( \frac{\lVert \mathbf{1}-\tilde{x}'\rVert_1}{ |\mathcal{P}| - N } \right)$ is a constant factor across $i\in\mathcal{P}$ in (\ref{e:transform}), selecting the smallest $N$ values of $1-\tilde{l}'_i$ is equivalent to selecting the smallest $N$ values of $1-\tilde{l}_i = \left(1-\tilde{x}_i\right) \left( 1 - \tilde{\pi}_i\right)$, which is equivalent to $x^{\mathrm{INT}}\left(\tilde{x}\right)$ selecting the largest $N$ values of $\tilde{l}_i = \tilde{x}_i + \left( 1 - \tilde{x}_i \right)\tilde{\pi}_i$. This completes the proof.
\qed
\end{pf}

During the convergence to the best known integer solution, a decelerator slows the speed of convergence. To determine $\alpha (t)$, we define the root-mean-square deviation (RMSD) of $\tilde{x}^{(t)}$ with respect to 0.5, {\it i.e.},
\begin{equation}
    \mathrm{RMSD}(\tilde{x}^{(t)}) = \sqrt{\frac{\sum_{i\in\mathcal{P}} \left(\tilde{x}^{(t)}_i-0.5\right)^2}{|\mathcal{P}|}},
\end{equation}
which ranges from 0 at $\ddot{x}$ to 0.5 at an integer solution $x^{\mathrm{INT}}$. Thus, as $\tilde{x}^{(t)}$ converges to the best known integer solution $x^{\mathrm{BKI}}$, the RMSD converges to 0.5.
The smoothing constant $\alpha (t)$ is determined by decelerator (DEC) which slows $\alpha = \frac{1}{2}$ (at $\ddot{x}$) to $\alpha = \frac{1}{1+e^2}$ (at $x^{\mathrm{BKI}}$):
\begin{equation}
    \label{e:decRR}
    \alpha (t) = \mathrm{DEC} (t-1) = \frac{1}{1+e^{4*\mathrm{RMSD}(\tilde{x}^{(t-1)})}}.
\end{equation}


\begin{algorithm}
\caption{Adaptive Randomized Rounding}
\label{a:advRR}
\begin{algorithmic}[1]
\STATE {\bf Initialize} $\tilde{x}^{(0)} \leftarrow \ddot{x}=\left(\ddot{x}_i= 0.5\mbox{ for }i\in\mathcal{P}\right)$; $x^{\mathrm{BKI}}\leftarrow x^{\mathrm{INT}}\left(\ddot{x}\right)$; 
$t \leftarrow 1$;
$n^\mathrm{local} \leftarrow 0$
\WHILE{all termination criteria are not met}
    \STATE Adjust smoothing constant $\alpha (t)$ by Decelerator (\ref{e:decRR})\label{l:alpha}
    \STATE Move $\tilde{x}^{(t)}$ by exponential smoothing (\ref{eqn:expSmoothNew}) with the smoothing constant $\alpha (t)$
    \STATE Find $x^{\mathrm{INT}(t)}=x^{\mathrm{INT}}\left(\tilde{x}^{(t)}\right)$
    \IF{$\mathrm{DEMAND}\left(x^{\mathrm{INT}(t)}\right)>\mathrm{DEMAND}\left(x^{\mathrm{BKI}}\right)$}
        \STATE $x^{\mathrm{BKI}} \leftarrow x^{\mathrm{INT}(t)}$ 
        \STATE $n^\mathrm{local} \leftarrow 0$
    \ELSIF{$x^{\mathrm{INT}(t)}=x^{\mathrm{BKI}}$}
        \STATE $n^\mathrm{local} \leftarrow n^\mathrm{local}+1$
        \STATE $p^\mathrm{R} \leftarrow \mathrm{min}(\frac{n^\mathrm{local}}{20}, 1) *\mathrm{RMSD}(\tilde{x}^{(t)})$
        \STATE Reset $\tilde{x}^{(t)} \leftarrow \ddot{x}$ and $n^\mathrm{local} \leftarrow 0$ with probability of $p^\mathrm{R}$
    \ELSE 
        \STATE $n^\mathrm{local} \leftarrow 0$
    \ENDIF
    \STATE $t \leftarrow t+1$
\ENDWHILE
\STATE Return the best known integer solution $x^{\mathrm{BKI}}$
\end{algorithmic}
\end{algorithm}

Pseudo-code for this ARR procedure is given in Algorithm~\ref{a:advRR}. The number of consecutive trials in which the same solution $x^{\mathrm{BKI}}$ is found is denoted by $n^{\mathrm{local}}$. As the same solution $x^{\mathrm{BKI}}$ is found in the run of trials ({\it i.e.}, $n^{\mathrm{local}}$ increases), the probability that the seed will be reset to $\ddot{x}$ increases. In particular, the probability of reset $p^\mathrm{R}$ reaches $\mathrm{RMSD}\left(\tilde{x}^{(t)}\right)$ after a run of 20 trials that find the same $x^{\mathrm{BKI}}$.
A computational time limit is set as the termination criterion for the full procedure. In the computational experiments on medium-scale instances, finding the global optimal solution was set as an additional termination criterion, as described in the next section.

\section{Numerical Experiments}\label{sec:experiments}
\subsection{{Experimental Design}}
Two numerical experiments, on medium- and large-scale instances, respectively, were conducted to evaluate the performance of the proposed algorithms (the neighborhood search (NS) method and the adaptive randomized rounding (ARR) procedure) and to compare the P\&R FLPs as modeled using the MNL and NL demand functions. In the experiments, the logsum parameter $\lambda$ of the NL demand function was set to 0.5 except when the $\lambda$ sensitivity analysis was performed. {As references for the performance comparison, three other algorithms were implemented on the large-scale instances: the conventional ARR procedure introduced by \citet{chopra2019}, denoted as Con-ARR; an ARR with an accelerator (ACC) (instead of a decelerator (DEC)), denoted as Acc-ARR; and a genetic algorithm (GA).}

Acc-ARR is the adaptive randomized rounding procedure of Algorithm~\ref{a:advRR} with the smoothing constant $\alpha(t)$ in Line~\ref{l:alpha} determined by accelerator ACC:
\begin{equation}
    \label{e:accRR}
    \alpha (t) = \mathrm{ACC} (t-1) = \frac{1}{1+e^{-8(\mathrm{RMSD}(\tilde{x}^{(t-1)})-0.5)}}.
\end{equation}
The smoothing constant increases from $\frac{1}{1+e^4}$ at $\ddot{x}$ to $\frac{1}{2}$ at $x^{\mathrm{INT}}$.
Details of Con-ARR and GA are presented in \ref{s:Con-ARR} and \ref{sec:ga}, respectively.

In both experiments, a P\&R FLP was solved for a randomly generated metropolitan area that consisted of a number of central business districts (CBDs) and neighborhoods (residential areas). 
The number of CBDs was a random number between 3 and 5, and the CBDs were randomly located on a radius between 1 and 2 distance units from the center of the metropolitan area. The number of neighborhoods was a random number between 5 and 10, and the center of each neighborhood was randomly located on a radius between 6 and 10 distance units from the center of the metropolitan area. Note that the distance unit is arbitrary; it could be taken as a kilometer or a mile, for example. We generated 1,000 random instances of the P\&R FLP and executed the proposed algorithms for a Monte Carlo statistical analysis.

The origin of each trip was located in a neighborhood and was uniformly distributed in an area $[-1,1]\times[-1,1]$ around the neighborhood center; thus, a neighborhood is a cluster of origins. The destination of each trip was assumed to be a CBD. Thus, multiple origins were paired with the same destination. P\&R facilities had to be located between the neighborhoods and the CBDs, and the candidate P\&R locations for a random instance were located on a radius between 5 and 7 distance units from the center of the metropolitan area. 

\begin{figure}[htb]
    \centering
    \includegraphics[width=\textwidth]{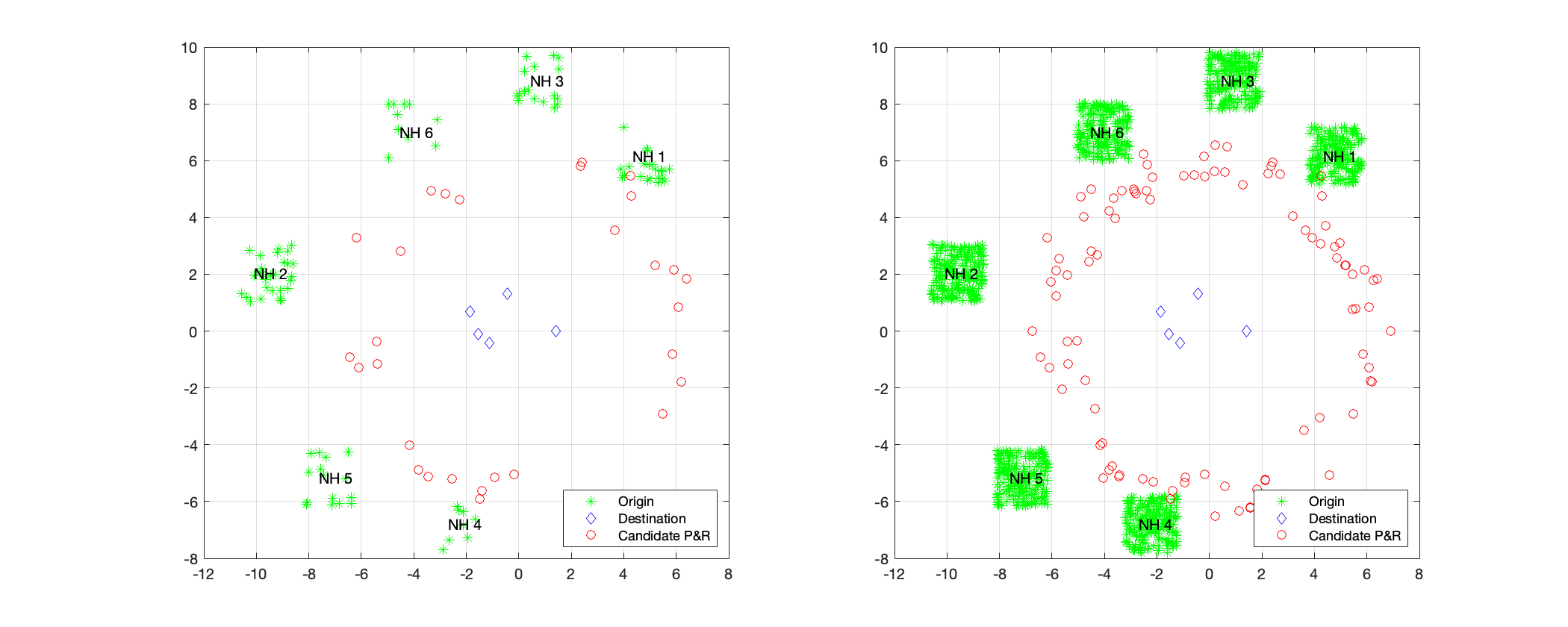}
    \caption{Illustrative Example of Origin, Destination, and Candidate P\&R Locations for Numerical Experiments}
    \label{f:od2}
\end{figure}

In Experiment~1, the P\&R FLP was solved on 1,000 medium-scale instances. The number of O-D pairs was 40, and eight P\&R locations were selected from among 30 candidate locations. In Experiment~2, the P\&R FLP was solved on 100 large-scale instances, in which the number of O-D pairs was 1,000, and 35 P\&R locations were selected from among 100 candidate locations. The two experiments differed in problem size ({\it i.e.}, the number of O-D pairs and candidate P\&R locations), but their geometric structures were the same. Figure~\ref{f:od2} shows an illustrative example of six neighborhoods and five CBDs for Experiment~1 on the left and Experiment~2 on the right.

In the experiments, the proposed algorithms were terminated when a time limit was reached. 
The time limits for Experiment~1 and Experiment~2 were set to 1 min and 30 min, respectively, based on preliminary trial runs.
Table~\ref{t:experiment} summarizes the conditions for each experiment. Because of the large search space for Experiment~2, a brute-force search was not applicable. 

\begin{table}[hbt!]
    \centering
    \caption{Summary of Numerical Experiment Details}
    \begin{tabular}{c|c|c}
    \hline
    & Experiment 1 & Experiment 2 \\ \hline
    Number of O-D pairs & 40 & 1,000 \\
    Number of P\&R candidates & 30 & 100 \\
    Number of selected P\&R locations & 8 & 35 \\
    Brute-Force search & Yes & No \\
    Time limit of heuristic algorithms & 1 min & 30 min \\
    Number of instances & 1,000 & 100 \\
    \hline
    \end{tabular}
    \label{t:experiment}
\end{table}

The experiments were conducted using MATLAB R2020b on a computer with an Intel Core i9-9900T CPU and 16GB of RAM. 
For the experiments, the negative of the travel distance (Euclidean distance) was used as the observed utility for simplicity of calculation. Thus, the longer the distance traveled, the less attractive the alternative. Thus, the utility values for the car option and the P\&R option are, respectively,
\begin{align}
    V_j^\mathrm{c} &= - TD^\mathrm{c}_j = -\sqrt{(X^\mathrm{d}_j-X^\mathrm{o}_j)^2+(Y^\mathrm{d}_j-Y^\mathrm{o}_j)^2}, ~\forall j \in \mathcal{T} \\
    V^\mathrm{p}_{ij} &= - TD^\mathrm{p}_{ij} = -\sqrt{(X^\mathrm{p}_i-X^\mathrm{o}_j)^2+(Y^\mathrm{p}_i-Y^\mathrm{o}_j)^2} -\sqrt{(X^\mathrm{d}_j-X^\mathrm{p}_i)^2+(Y^\mathrm{d}_j-Y^\mathrm{p}_i)^2}, ~\forall i \in \mathcal{P}, ~\forall j \in \mathcal{T},
\end{align}
where ($X^\mathrm{o}_j$, $Y^\mathrm{o}_j$); and ($X^\mathrm{d}_j$, $Y^\mathrm{d}_j$) are the X- and Y-coordinates of the origin and destination of O-D trip $j$, and ($X^\mathrm{p}_i$, $Y^\mathrm{p}_i$) are those for P\&R candidate location $i$. Note that O-D trip $j$ comprises $\{(X^\mathrm{o}_j, Y^\mathrm{o}_j), (X^\mathrm{d}_j, Y^\mathrm{d}_j)\}$, and the utility of a private car is always greater than that of any P\&R option. The number of travelers on each O-D trip was set to one for simplicity; {\it i.e.},
$$R_j = 1, ~\forall j \in \mathcal{T}.$$

\begin{figure}[htb]
    \centering
    \includegraphics[width=0.7\textwidth]{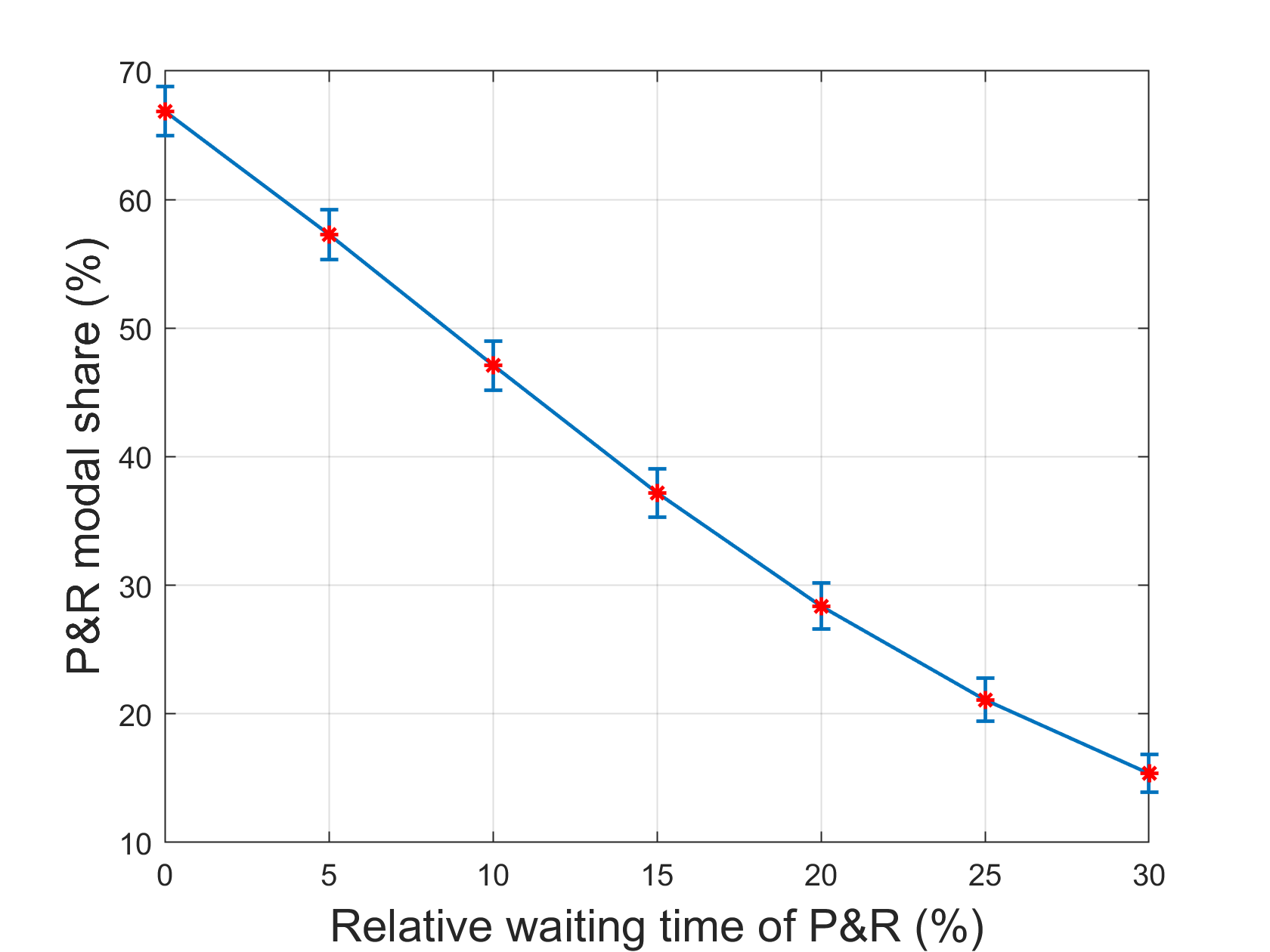}
    \caption{Impact of Waiting Time on the Modal Share of P\&R}
    \label{f:saWaiting}
\end{figure}

{As given in (\ref{e:Vcar}) and (\ref{e:Vpr}), the utility (and therefore the modal share) of the P\&R facility may be affected by other factors, such as the wait time at a P\&R facility and the frequency of the transit service. For example, Fig.~\ref{f:saWaiting} shows the impact of wait time on the modal share of the P\&R option in Experiment~2. The wait time may be added to the utility of the P\&R option as follows:
$$V^\mathrm{p}_{ij} = - TD^\mathrm{p}_{ij}*(1+\Delta W),$$ 
where $\Delta W$ is the relative wait time at the P\&R facility and can range from 0 to 0.3. As expected, the modal share of the P\&R option decreases as the wait time increases because the P\&R utility decreases. Likewise, other factors included in the P\&R utility may influence the decision making of potential P\&R users. However, we did not observe any unexpected effects of other factors in preliminary experiments, and therefore we considered only the travel distance in the experiments. Thus, the negative of the travel distance was used as the utility for both the private car mode and the P\&R mode. For a MATLAB code and the data for the numerical experiments, the readers may refer to the data repository provided in \citet{mydata}.}

\subsection{{Comparison of Algorithm Performance}}
To evaluate and compare the performance of the proposed algorithms, we performed the two experiments described in the previous section. For Experiment~1, the number of feasible solutions was $\binom{30}{8} = 5,852,925$, and a brute-force search was computationally tractable for this medium-scale problem. We used the exact optimal solution as the benchmark to evaluate the quality of the solution produced by the proposed algorithms. In contrast, for Experiment~2, the number of feasible solutions was $\binom{100}{35}$, which is approximately $1.87\times e+20$ times greater than that for Experiment~1. As the number of candidate P\&R locations increases, the search space increases exponentially, and a brute-force search cannot solve the large-scale instances to exact optimality within a reasonable time. 
Therefore, instead of using the exact optimal solution, for Experiment~2 we used the best solution found as the benchmark. 

Table~\ref{t:perf1} shows the computational performance of the proposed algorithms and the comparison algorithms in Experiment~1. All algorithms successfully found the exact optimal solutions to all of the medium-scale instances (``Optimality"). {NS initially converged to sub-optima on seven of the 1,000 instances. Within the computational time limit, the subsequent NS was performed with a random initial solution $x^{\mathrm{INT}}\left(\ddot{x}\right)$. All seven sub-optima were improved to the exact optimal solutions after two or more repetitions of the NS with the random initial solutions.} The brute-force (BF) search provided the exact optimal solutions as the benchmark for these medium-scale instances, and the termination criterion for all of the algorithms was the exact optimality of the solution. The ``Mean Time," ``Std Time," ``Max Time," and ``Mean Iterations" rows show the average, standard deviation, and maximum computation time and the average number of iterations on the 1,000 instances until exact optimality was reached. For NS, the number of iterations is the number of edges traversed; for ARR, Acc-ARR, and Con-ARR, it is the number of trials; and for GA, it is the number of generations. NS and ARR were much faster than the brute-force search. They also outperformed Acc-ARR, Con-ARR, and GA in computational time.

\begin{table}[hbt!]
    \centering
    \caption{Comparison of the {\bf Proposed Algorithms} in Experiment 1 on One Thousand Instances}
    \begin{tabular}{c|c|c|c|c|c|c}
    \hline
    & BF & {\bf NS} & {\bf ARR} & Acc-ARR & Con-ARR & GA \\ \hline
    Optimality & 100\% & {\bf 100\%} & {\bf 100\%} & 100\% & 100\% & 100\% \\
    Mean Time & 158 s & {\bf 0.04 s} & {\bf 0.05 s} & 0.19 s & 0.12 s & 0.31 s \\
    Std Time & 1.31 s & {\bf 0.01 s} & {\bf 0.17 s} & 0.55 s & 0.21 s & 0.46 s \\
    Max Time & 165 s & {\bf 0.10 s} & {\bf 3.6 s} & 15 s & 5.9 s & 7.5 s \\
    Mean Iterations & N/A & {\bf 954} & {\bf 1,012} & 3,742 & 2,190 & 522 \\
    \hline
    \end{tabular}
    \label{t:perf1}
\end{table}

Table~\ref{t:perf2} shows the computational performance of the proposed algorithms and the comparison algorithms in Experiment~2 for the time limit of 30 min. As explained above, the brute-force search cannot solve the large-scale instances, and thus there is no benchmark solution for exact optimality. The benchmark solution for each instance was the best solution found (BSF) by all five algorithms. 
``Number of instances with the BSF" shows that the proposed algorithms reached the BSF on all 100 instances along with Acc-ARR and Con-ARR, whereas GA did not reach any BSF. In fact, all but GA together reached the same BSF on every instance. ``Mean ratio of demand to the BSF" shows the ratio of the best demand found by the algorithm to the demand according to the BSF; the average gap between the GA solution and the BSF was 3.5\%. The computational time for each algorithm to reach its best solution is denoted in the table by $t^\textrm{f}$. The two proposed algorithms reached the BSF much earlier than the computational time limit of 30 min, and each was faster than Acc-ARR, Con-ARR, and GA. In particular, they reached the BSF much faster than GA reached its own best solutions, which were inferior to the BSF. ``Mean iterations until time limit" shows the average number of iterations of the algorithm within the computational time limit of 30 min. NS found the BSF on each of the large-scale instances without repetition ({\it i.e.}, beginning with random solution $x^{\mathrm{INT}}\left(\ddot{x}\right)$ only once). As the number of iterations within the time limit is nearly equal across ARR, Acc-ARR, and Con-ARR, the average number of iterations to reach the BSF is proportional to the mean $t^\textrm{f}$ for each of the three randomized rounding procedures.

\begin{table}[hbt!]
    \centering
    \caption{Comparison of the {\bf Proposed Algorithms} in Experiment 2 on 100 Instances (Time Limit = 30 minutes)}
    \begin{tabular}{c|c|c|c|c|c}
    \hline
    & {\bf NS} & {\bf ARR} & Acc-ARR & Con-ARR & GA \\ \hline
    Number of instances with the BSF & {\bf 100} & {\bf 100} & 100 & 100 & 0 \\
    Mean ratio of demand to the BSF & {\bf 1} & {\bf 1} & 1 & 1 & 0.965 \\
    Mean $t^\textrm{f}$ 
    & {\bf 29 s} & {\bf 70 s} & 271 s & 149 s & 1,037 s \\
    Std $t^\textrm{f}$ 
    & {\bf 10 s} & {\bf 64 s} & 219 s & 169 s & 525 s \\
    Max $t^\textrm{f}$ 
    & {\bf 60 s} & {\bf 316 s} & 1,284 s & 845 s & 1,793 s \\
    Mean iterations until time limit & {\bf 16,088} & {\bf 552,695} & 552,736 & 552,206 & 31,581 \\
    \hline
    \end{tabular}
    \label{t:perf2}
\end{table}

\begin{figure}[htb]
    \centering
    \includegraphics[width=0.7\textwidth]{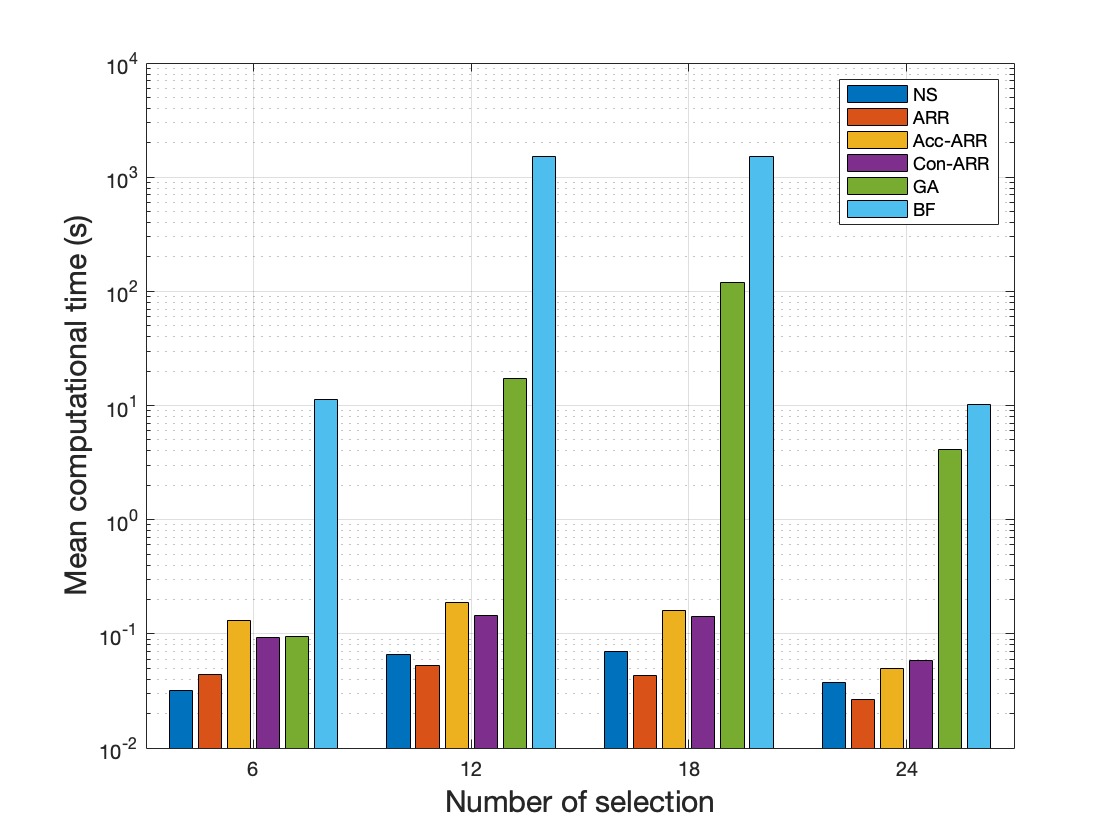}
    \caption{Impact of $N$ on the Average Convergence Speed}
    \label{f:saN}
\end{figure}

{
We performed a sensitivity analysis to understand the effect of various numbers of selected P\&R locations ($N$). We chose 100 random instances of Experiment~1 with 40 O-D and 30 candidate P\&R locations and applied the algorithms to these instances to select the optimal $N$ locations that maximize demand, for $N$ values of 6,12,18, and 24. All of the algorithms except GA found the exact optimal solutions for all $N$ values. Figure~\ref{f:saN} shows the average computational time taken by the algorithms (corresponding to ``Mean Time" in Table~\ref{t:perf1}) for the four values of $N$. The Y axis shows the computational time on a logarithmic scale. As expected, BF took significantly more computational time than did the other algorithms. The computational times of BF for $N$ values of 6 and 24 are nearly equal because $\binom{30}{6} = \binom{30}{24}$; similarly, the computational times for $N$ values of 12 and 18 are nearly equal. The computational time of BF for an $N$ value of 12 or 18 is much greater than that for an $N$ of 6 or 24 because $\binom{30}{12} > \binom{30}{6}$. The figure also shows that NS and ARR were faster than the other algorithms for the various value of $N$, ARR begin slightly faster than NS for the $N$ values of 12, 18, and 24.
}

{
Now, we compare the performance of the proposed algorithms with that of the MILP formulation introduced by \citet{aros2013p} on 100 medium-scale P\&R FLPs (from Experiment~1) under the MNL demand function.
The P\&R FLPs select 15 P\&R locations from among 30 candidate locations, on which BF is intractable because of the large search space.
The MILP (\ref{e:obj2})-(\ref{e:var2}) was solved by Gurobi Optimizer version 9.0.3, and NS and ARR solved (\ref{e:brief}) with $\lambda=1$. Table~\ref{t:perf3} shows the comparison of the performance of the algorithms. The proposed algorithms solved the medium-scale MNL instances to exact optimality 10,000 times faster than the MILP ($\mathrm{MILP}: \mathrm{NS} = 295\mbox{s}: 0.02\mbox{s}$) where the computational time for MILP was the system time excluding the time used for construction and reading of the model. Moreover, the proposed algorithms solved the large-scale MNL models (from Experiment~2) within two minutes, similar to the performance reported in Table~\ref{t:perf2}, whereas the MILP is intractable on the large-scale MNL models generating the error message of Out-Of-Memory.
}

\begin{table}[hbt!]
    \centering
    \caption{Comparison of Commercial Solver (MILP) and the {\bf Proposed Algorithms} in Experiment~1 with MNL model}
    \begin{tabular}{c|c|c|c}
    \hline
    & MILP & {\bf NS} & {\bf ARR}  \\ \hline
    Optimality & 100\% & {\bf 100\%} & {\bf 100\%} \\
    Mean computational time & 295 s & {\bf 0.02 s} & {\bf 0.06 s} \\
    Std computational time & 484 s & {\bf 0.01 s} & {\bf 0.05 s} \\
    Max computational time & 2,639 s & {\bf 0.05 s} & {\bf 0.2 s} \\
    \hline
    \end{tabular}
    \label{t:perf3}
\end{table}

\subsection{Comparison of Logit Models}\label{s:comparelogits}
The second purpose for conducting the numerical experiments was to analyze the dependency of the optimal P\&R locations on the demand modeling. As stated previously, the MNL model and the NL model are based on different assumptions regarding the correlation between P\&R locations. {Specifically, the MNL model assumes the logsum parameter $\lambda$ to be 1, whereas the NL model estimates $\lambda$ by using maximum likelihood estimation \citep{koppelman1998alternative}.} As no {\it a priori} knowledge of $\lambda$ for the P\&R FLP was available, the experiments reported in the previous section arbitrarily assumed $\lambda=0.5$. Because $\lambda$ is an important parameter indicating the underlying correlation between P\&R locations \citep{koppelman2006self}, we conducted a sensitivity analysis of $\lambda$ on the selected P\&R locations and the P\&R demand. ARR was used to solve the same 100 instances sampled from the 1,000 instances of Experiment~1 with $\lambda$ values ranging from 0.1 to 1 in steps of 0.1.

\begin{figure}[htb]
\centering
\includegraphics[width=0.7\textwidth]{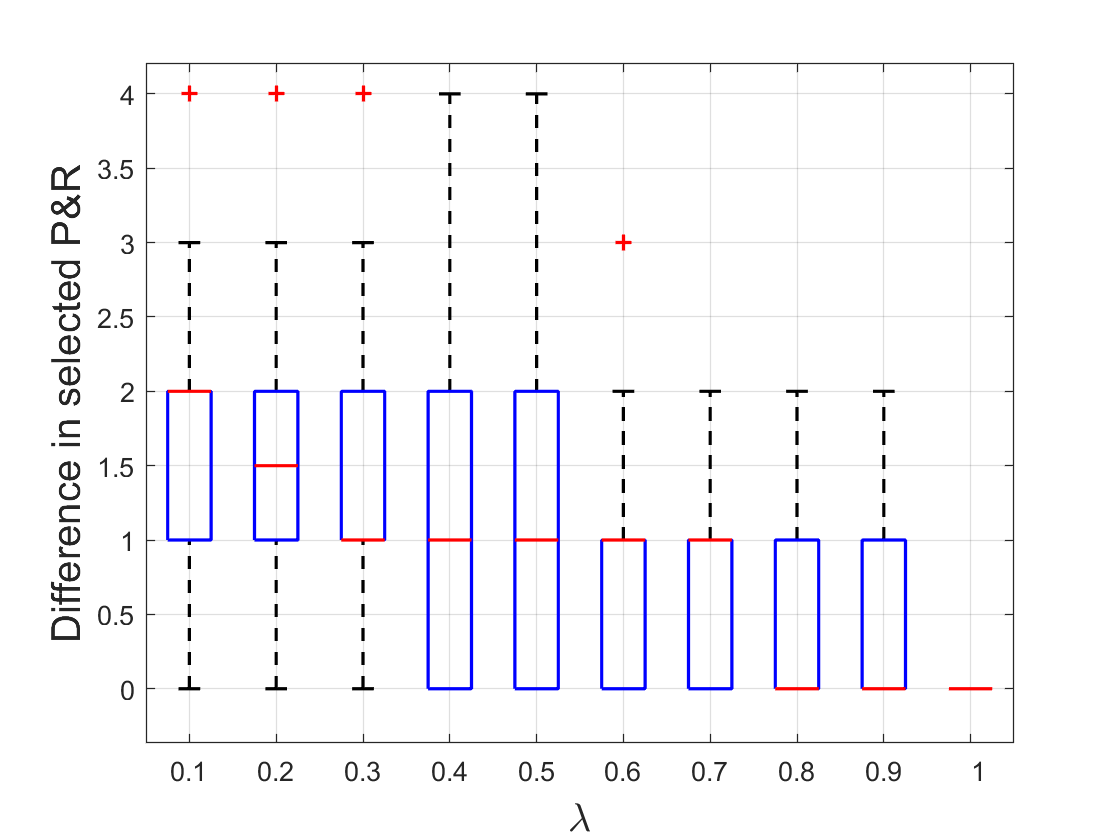}
\caption{Box Plot of Difference in 8 P\&Rs Selected from 30 Candidates with Respect to MNL Model}
\label{f:saDiff}
\end{figure}

Figure~\ref{f:optimal} shows a comparison of the probabilities associated with choosing each option (the private car option and one of the P\&R options) and the optimal combination of P\&R locations for an instance of Experiment~1. The $N$ P\&R locations selected as optimal are marked with an asterisk. The figure shows that for the same P\&R FLP, the two models gave not only different values for the P\&R demand but also different combinations of the $N$ P\&R locations. As can be seen, P\&R locations 5, 7, 15, 18, 22, and 28 were selected by both models, but the models differed in their selections for two P\&R locations. Figure~\ref{f:saDiff} shows the differences in the eight P\&R locations selected from among 30 candidates for various value of $\lambda$ with respect to the MNL model ({\it i.e.}, for $\lambda = 1$). It was found that the difference increases as the correlation between P\&R locations increases ({\it i.e.}, as $\lambda$ decreases). Note that the red bars in the figure are median values 
and the red crosses are outliers.

\begin{figure}[htb]
\centering
\includegraphics[width=0.7\textwidth]{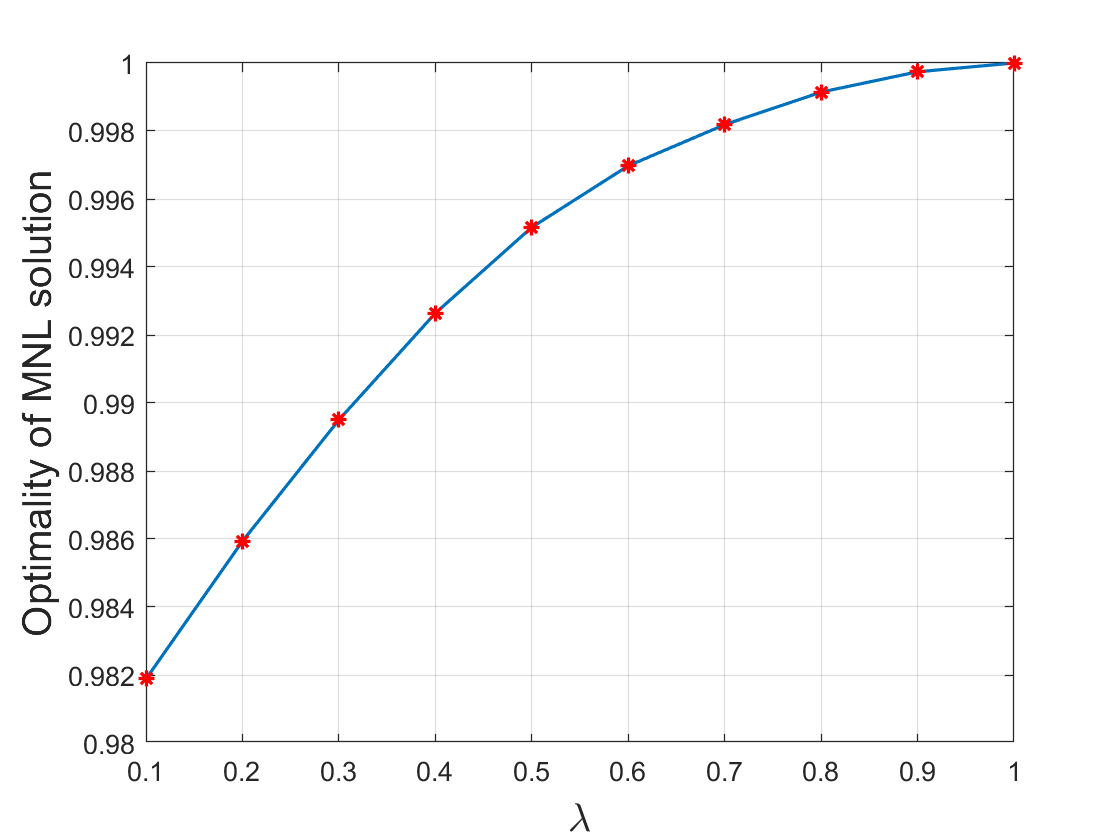}
\caption{Optimality of MNL Solution with Respect to NL Model with $\lambda$}
\label{f:mnlOpt}
\end{figure}

Because of the IIA assumption of the MNL model, the probability that a private car is chosen is lower in the MNL model than when the NL model is used. That is, the P\&R demand is overestimated by the MNL model, and a solution of the P\&R FLP with the MNL model will be sub-optimal, considering the correlation between P\&R locations. {Fig.~\ref{f:mnlOpt} shows the optimality of the MNL solution for the NL model with various values of $\lambda$. Let $x^{\mathrm{NL}(\lambda)}$ be the optimal solution to the NL model with $\lambda$. Then, the optimality of the MNL solution $x^{\mathrm{NL}(1)}$ is
\begin{eqnarray*}
\frac{\mathrm{DEMAND}_{\lambda}\left( x^{\mathrm{NL}(1)}\right)}{\mathrm{DEMAND}_{\lambda}\left( x^{\mathrm{NL}(\lambda)}\right)}.
\end{eqnarray*}
For example, suppose that there is a correlation between P\&R locations and that the corresponding value of $\lambda$ is 0.2. For $\lambda = 0.2$, the optimality gap is approximately 1.4\%. For $\lambda < 0.33$, the gap is greater than 1\%. The greater the correlation between P\&R locations, the greater the optimality gap of the MNL solution.}

\section{Conclusion}\label{sec:conclusion}
{
This study proposed two heuristic algorithms that can efficiently solve a park-and-ride facility location problem under the nested logit demand function. 
Previously the park-and-ride facility location problem was solved under the multinomial logit demand function in mixed-integer linear programming (MILP) formulations. The multinomial logit demand function does not capture the nested structure of the private car and park-and-ride options, and the MILP formulations can solve only the medium-scale park-and-ride facility location problem. Our heuristics of neighborhood search and adaptive randomized rounding were able to solve the large-scale instances extremely rapidly under the nested logit demand function, which captures the nested structure of the alternatives. In particular, numerical experiments verified that the proposed heuristics solve the medium-scale problem under the multinomial logit demand function to exact optimality much faster than the MILP formulations.
Moreover, the proposed algorithms can identify the best solutions in minutes for large-scale problems that are intractable under the MILP approach. 
}

{
The two heuristic algorithms proposed by this paper can solve a large-scale park-and-ride facility location problem to select a large number of new transit facilities for a new mode of transportation such as urban air mobility.
In particular, our adaptive randomized rounding procedure is not only for the polyhedral structure of a specific formulation, and can address more complicated facility location problems as well.
For example, a vertiport location problem may select the optimal pairs of vertiports for the commuters who use urban air mobility from one vertiport to another. 
The adaptive randomized rounding procedure is also expected to efficiently solve many other large-scale vertiport location problems. 
}



\section*{Acknowledgments}
The authors gratefully thank anonymous referees for helpful comments on the earlier version of the paper. We would like to thank Editage (www.editage.co.kr) for English language editing. 
This work is supported by the Korea Agency for Infrastructure Technology Advancement grant funded by the Ministry of Land, Infrastructure and Transport (Grant RS-2022-00143965). 

\appendix

\begin{figure}[htb]
    \centering
    \includegraphics[width=0.9\textwidth]{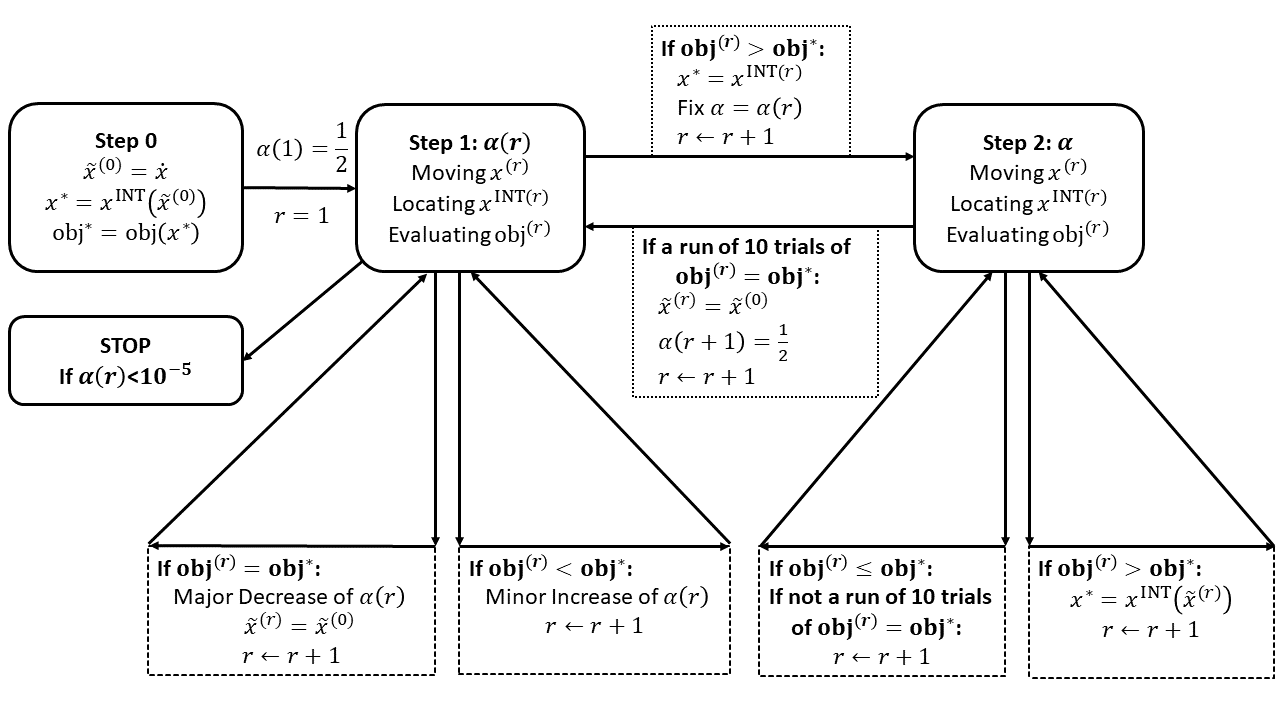}
    \caption{Flow Chart of 2 Step Procedure (where $\mathrm{obj}^{(r)}=\mathrm{obj} ( x^{\mathrm{INT}(r)} )$ and $\mathrm{obj}^* = \mathrm{obj}\left( x^{\mathrm{BKI}}\right)$) where $\mathrm{obj}=\mathrm{DEMAND}$}
    \label{f:twoStep}
\end{figure}

\section{Conventional adaptive randomized rounding}\label{s:Con-ARR}
Setting $\tilde{x}^{(0)}=\dot{x}\in\mathrm{REL}(N)$, the following $\tilde{x}^{(r)}$, $r>0$, will belong to $\mathrm{REL}(N)$ because the polyhedron contains the line segment $\left[ x^{(r-1)}, x^{\mathrm{BKI}}\right]$.
Moving $\tilde{x}^{(r)}$ in $\mathrm{REL}(N)$, our randomized rounding method performs the following 2 step procedure (see Fig.~\ref{f:twoStep}):

\paragraph{Step 0. Initialization}
With an initial distribution $\tilde{x}^{(0)}$ and an observed $x^{\mathrm{BKI}} = x^{\mathrm{INT}(0)} = x^{\mathrm{INT}} \left(\tilde{x}^{(0)}\right)$, 
Step 0 sets the initial smoothing constant $\alpha (1) = 1/2$ and moves on to Step 1.

\paragraph{Step 1. Adjusting smoothing constant}
The first step moves the last distribution $\tilde{x}^{(r-1)}$ to $\tilde{x}^{(r)}$ on the line segment $\left[ \tilde{x}^{(0)}, x^{\mathrm{BKI}} \right]$ by exponential smoothing (\ref{eqn:expSmoothNew}), adjusting the smoothing constant $\alpha$, and observes a random combination $x^{\mathrm{INT}(r)}$ of $N$ facilities using distribution $\tilde{x}^{(r)}$. 
\begin{itemize}
    \item Stopping criterion of the whole procedure: If $\alpha (r)$ is less than a certain tolerance (e.g., tolerance = $10^{-5}$ equal to the integrality tolerance of \citet{gurobi}), the whole procedure stops. 
    \item Criterion to go to Step 2: If  $\mathrm{obj}\left(x^{\mathrm{INT}(r)}\right)>\mathrm{obj}\left(x^{\mathrm{BKI}}\right)$, the procedure updates the best known combination $x^{\mathrm{BKI}} = x^{\mathrm{INT}(r)}$, where $\mathrm{obj}=\mathrm{DEMAND}$. Assuming that the speed $\alpha (r)$ is appropriate, fix the present smoothing constant $\alpha = \alpha (r)$, set $r\leftarrow r+1$ and go to Step 2.
    \item Major decrease of $\alpha$: If $\mathrm{obj}\left(x^{\mathrm{INT}(r)}\right)=\mathrm{obj}\left(x^{\mathrm{BKI}}\right)$, assume that $\tilde{x}^{(r)}$ is arriving too close at $x^{\mathrm{BKI}}$ or too fast, and slow down the speed of the move by major decrease of the smoothing constant $\alpha (r+1) = \alpha (r) / 2$. Reset the present distribution back to $\tilde{x}^{(0)}$, set $r\leftarrow r+1$ and repeat Step 1.
    \item Minor increase of $\alpha$: If $\mathrm{obj}\left(x^{\mathrm{INT}(r)}\right)<\mathrm{obj}\left(x^{\mathrm{BKI}}\right)$, assume that $\tilde{x}^{(r)}$ is moving slow and speed up the move a little bit by minor increase of the smoothing constant $\alpha (r+1) = 2^{\lfloor\log_2 \alpha (r)\rfloor} + \alpha (r) / 2$ such as $\alpha(r)=1/2+1/4$ to $\alpha (r+1) = 1/2+1/4+1/8$. Set $r\leftarrow r+1$ and repeat Step 1 (moving ahead without coming back to $\tilde{x}^{(r+1)}=\tilde{x}^{(0)}$).
\end{itemize}

\paragraph{Step 2. Fixed smoothing constant}
The second step moves $\tilde{x}^{(r)}$ from $\tilde{x}^{(r-1)}$ toward $x^{\mathrm{BKI}}$ by exponential smoothing (\ref{eqn:expSmoothNew}) on the line segment $\left[ \tilde{x}^{(r-1)},x^{\mathrm{BKI}}\right]$ without changing the smoothing constant $\alpha$ that is adjusted in Step 1. It observes $x^{\mathrm{INT}(r)}$ based on $\tilde{x}^{(r)}$.
\begin{itemize}
    \item Criterion to go back to Step 1: If $\mathrm{obj}\left( x^{\mathrm{INT}(r)}\right) = \mathrm{obj}\left( x^{\mathrm{BKI}}\right)$ for a run of say 10 iterations of trials, assume that $\tilde{x}^{(r)}$ is too close to $x^{\mathrm{BKI}}$ and reset $\tilde{x}^{(r)}=\tilde{x}^{(0)}$ initializing $\alpha (r+1)=1/2$. Set $t\leftarrow t+1$ and go back to Step 1.
    \item If $\mathrm{obj}\left( x^{\mathrm{INT}(r)}\right) > \mathrm{obj}\left( x^{\mathrm{BKI}}\right)$, the procedure updates $x^{\mathrm{BKI}}=x^{\mathrm{INT}(r)}$ with the new best value. Set $r\leftarrow r+1$ and repeat Step 2 from the present distribution $\tilde{x}^{(r-1)}$ toward the new best known integer solution $x^{\mathrm{BKI}}$.
    \item If $\mathrm{obj}\left( x^{\mathrm{INT}(r)}\right) \leq \mathrm{obj}\left( x^{\mathrm{BKI}}\right)$, set $r\leftarrow r+1$ and repeat Step 2 going on from the present $\tilde{x}^{(r-1)}$ toward the best known facility location $x^{\mathrm{BKI}}$.
\end{itemize}

\section{Genetic Algorithm}\label{sec:ga}
Genetic Algorithm (GA) is a population-based heuristic algorithm \citep{holland1992adaptation}. It is well-known and widely adopted to solve the FLP \citep{kratica2001solving}. A canonical GA is given in Algorithm \ref{a:ga} \citep{whitley1994genetic}. 

\begin{algorithm}
\caption{Canonical GA Algorithm}
\label{a:ga}
\begin{algorithmic}
\STATE (Initialization) Randomly generate an initial population
\WHILE{all termination criteria are not met}
\STATE (Evaluation) Evaluate fitness values of the population
\STATE (Selection) Randomly select chromosomes ({\it i.e.}, parents) based on fitness
\FOR{every two parents}
\STATE (Crossover) Mix the parent chromosomes to form a child chromosome
\ENDFOR
\STATE (Mutation) With probability of $p^\mathrm{m}$, randomly select genes for each child and change their values
\STATE (Next generation) chromosomes $\leftarrow$ children chromosomes
\ENDWHILE
\STATE Return the best objective value
\end{algorithmic}
\end{algorithm}

\paragraph{Phase 1: Initialization}
Each chromosome is a binary string of length $|\mathcal{P}|$, and $N$ genes ({\it i.e.}, elements of a chromosome) are randomly chosen and assigned to 1. Other genes are assigned to 0. Thus, a chromosome is a feasible set of $(x_1, x_2, \ldots, x_{|\mathcal{P}|})$. The population size ($\mathrm{pop}$) is a design parameter.

\paragraph{Phase 2: Evaluation}
The fitness value of a chromosome is equal to the objective value; {\it i.e.},
$$\mathrm{Fitness} = \sum_{i \in \mathcal{P}} \sum_{j \in \mathcal{T}} R_j p^\mathrm{p}_{ij},$$
but it can also be determined with the rank of objective value \citep{baker1985adaptive}. 

\paragraph{Phase 3: Selection}
The selection probability of a chromosome is proportional to its fitness, and the selection process is stochastic sampling with replacement, as in the roulette wheel method \citep{whitley1994genetic, karakativc2015survey, kim2019receding}. In a canonical GA, parents produce two children; however, in this paper, they produce one child. As a result, $2\times\mathrm{pop}$ parents are selected, and thus most chromosomes are selected multiple times as parents. 

\paragraph{Phase 4: Crossover}
A canonical GA utilizes 1-point crossover, which randomly selects a crossover point and swaps the fragmented genes between the parents. Some implementations of GA utilize 2-point crossover in which parents' middle genes between the crossover points are swapped to generate children \citep{andre2001improvement}. However, such a fragment-based crossover does not guarantee the feasibility of children chromosomes because (\ref{e:selection}) constrains the number of chosen P\&R to be equal to $N$.

\begin{figure*}
    \centering
    \includegraphics[width=0.9\textwidth]{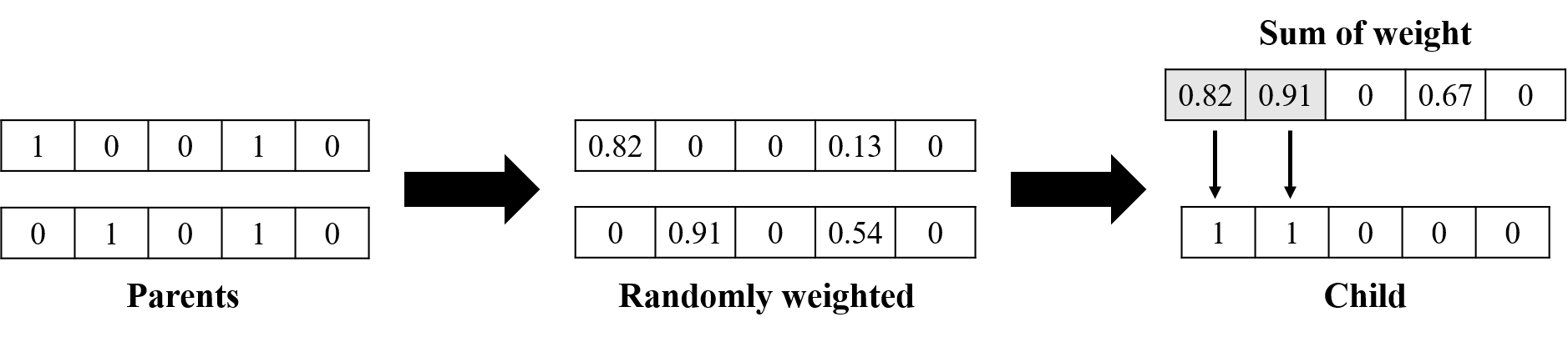}
    \caption{Weight-Based Crossover Operator}
    \label{f:cross}
\end{figure*}

Therefore, this paper proposes a weight-based crossover operator presented in Fig.~\ref{f:cross}. Each gene of parent chromosomes is randomly weighted between 0 and 1, and genes are added together to form a child chromosome. Then, the $N$ largest genes of the child chromosome are assigned to 1 and the other genes are assigned to 0. This crossover operator guarantees that the child chromosome is inherited only from its parents' chromosomes and the number of genes with value of 1 remains at $N$. Thus, the child chromosome is also a feasible solution of the P\&R FLP.

\paragraph{Phase 5: Mutation}
A canonical GA mutates each gene ({\it i.e.}, flips a bit from 0 to 1 or vice versa) with a probability of $p^\mathrm{m}$; however, such a mutation operator can make children chromosomes infeasible. Therefore, the mutation operator in this paper exchanges a gene with value of 1 with another gene with value of 0. Each child chromosome is subject to mutation with a probability of $p^\mathrm{m}$, which is a design parameter.

\begin{figure}[h]
    \centering
    \includegraphics[width=0.7\textwidth]{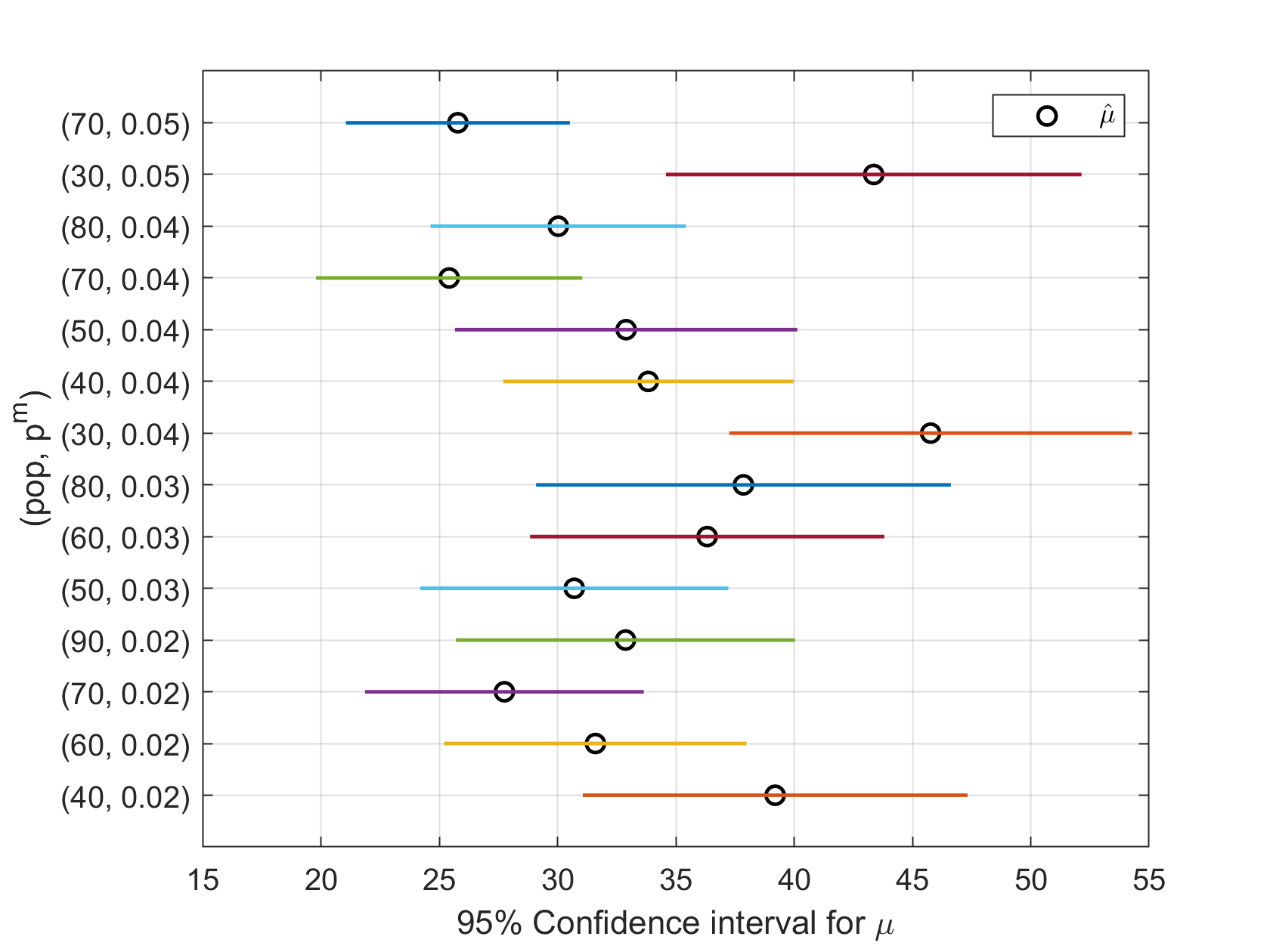}
    \caption{$\hat{\mu}$ and 95\% Confidence Interval for $\mu$ in GA Design Parameter Tuning}
    \label{f:ci}
\end{figure}

Design parameters of GA influence the performance, specifically in terms of optimality and convergence time \citep{angelova2011tuning}. A canonical GA has design parameters such as population size, crossover probability, and mutation probability. This paper sets the crossover probability at 1, meaning all children chromosomes are inherited from both parents. However, the other parameters such as population size and mutation probability should be determined by design parameter tuning.

Prior to the numerical experiments, Monte-Carlo experiments were conducted to tune the population size ($\mathrm{pop}$) and mutation probability ($p^\mathrm{m}$). Based on previous experience, $\mathrm{pop}$ ranges from 30 to 90 in steps of 10, and $p^\mathrm{m}$ ranges from 1\% to 5\% in steps of 1\%. 100 instances of Experiment 1 are solved by GA for all pairs of $(\mathrm{pop},~p^\mathrm{m})$. The exact optimal solution to instances is found using a brute-force search, and the optimality of a GA solution is the ratio of the GA solution to the brute-force solution. Because the global optimal solution ({\it i.e.}, the brute-force solution) is known, the computational times at which the GA finds the global optimal solution ($t^\mathrm{f}$) are recorded. The true mean ({\it i.e.}, population mean) of $t^\textrm{f}$ is denoted as $\mu$, and its point estimate is denoted as $\hat{\mu}$. The time limit is set at 20 min based on sample runs.

\begin{figure}[htb]
    \centering
    \includegraphics[width=0.7\textwidth]{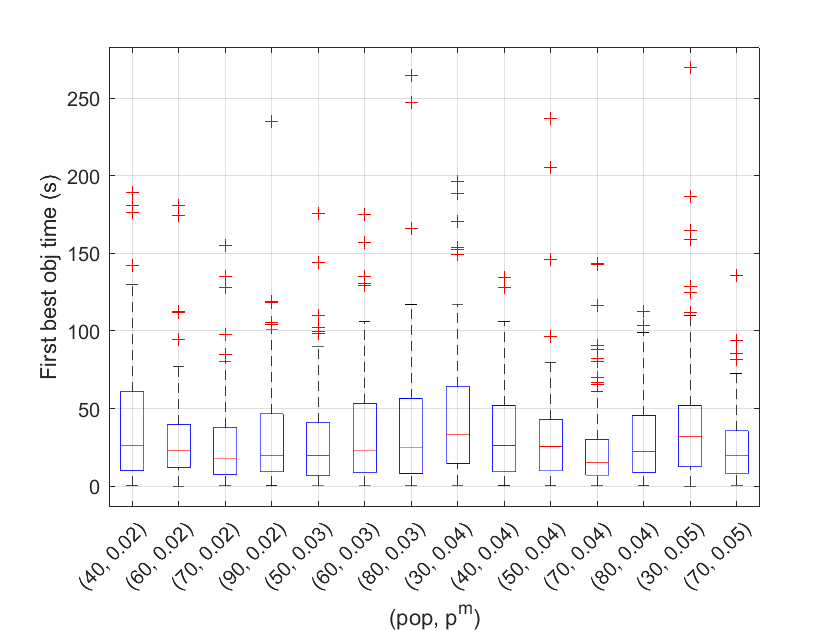}
    \caption{Box Plot of $t^\textrm{f}$ in GA Design Parameter Tuning}
    \label{f:box}
\end{figure}

Fig.~\ref{f:ci} shows $\hat{\mu}$ and 95\% confidence intervals for $\mu$, only for the pairs of $(\mathrm{pop},~p^\mathrm{m})$ with which the GA found the optimal solution to all instances within the time limit. Fig.~\ref{f:box} shows a box plot of $t^\textrm{f}$, also corresponding to $(\mathrm{pop},~p^\mathrm{m})$ pairs of 100\% optimality. Figs.~\ref{f:ci}--\ref{f:box} show that the (70, 0.04) pair has the fastest convergence time on average ({\it i.e.}, the smallest $\hat{\mu}$), as well as, in the worst case ({\it i.e.}, the smallest upper bound of $t^\textrm{f}$ excluding outliers). Hence, the parameter pair (70, 0.04) is used for the numerical experiments, but its 95\% confidence interval for $\mu$ overlaps with those for other $(\mathrm{pop},~p^\mathrm{m})$ pairs. For instance, (70, 0.05) can be another good parameter pair for GA.

\bibliography{flp}

\end{document}